\documentclass[11pt]{amsart}
\usepackage[margin=1in]{geometry}
\usepackage{amssymb,graphicx}
\usepackage{bm}
\usepackage[hidelinks]{hyperref}
\usepackage{cleveref}
\theoremstyle{plain}
\newtheorem{theorem}{Theorem}[section]
\newtheorem{lemma}[theorem]{Lemma}
\newtheorem{proposition}[theorem]{Proposition}

\theoremstyle{definition}
\newtheorem{fact}[theorem]{Fact}
\newtheorem{remark}[theorem]{Remark}

\newtheorem{example}[theorem]{Example}
\numberwithin{equation}{section}

\DeclareMathOperator{\expct}{\mathbb{E}}
\DeclareMathOperator{\tr}{tr}

\newcommand{\mat}[1]{\boldsymbol{#1}}
\newcommand{\vect}[1]{{\bf #1}}
\newcommand{\myfig}[2]{\includegraphics[width=#1\textwidth]{./figures/#2}}

\newcommand{\FG}[1]{Fig.~\ref{#1}}
\newcommand{\Thm}[1]{Theorem~\ref{#1}}
\newcommand{\SC}[1]{Sec.~\ref{#1}}
\newcommand{\expect}[1]{\expct\left[#1\right]}
\newcommand{\EP}{\expect{\proj}}
\newcommand{\EPxP}{\expect{\proj \otimes \proj}}
\newcommand{\proj}{\boldsymbol{P}}
\newcommand{\B}{\boldsymbol{\mathcal{B}}}
\newcommand{\C}{\boldsymbol{\mathcal{C}}}
\newcommand{\A}{\boldsymbol{\mathcal{A}}}
\newcommand{\T}{\boldsymbol{\mathcal{T}}}

\newcommand{\IxI}{\boldsymbol{\mathcal{I}}}
\newcommand{\I}{\boldsymbol{I}}
\newcommand{\norm}[1]{\lVert#1\rVert}

\newcommand{\reel}{\mathbb{R}}
\newcommand{\bx}{{\vect x}}

\newcommand{\bb}{{\vect b}}
\newcommand{\be}{{\vect e}}
\newcommand{\berr}{\mathbf{z}}

\newcommand{\mA}{{\mat A}}

\newcommand{\mE}{{\mat E}}

\newcommand{\mU}{{\mat U}}
\newcommand{\mT}{{\mat T}}
\newcommand{\mL}{{\mat L}}
\newcommand{\BB}{\overline{\varphi}_{\mathrm{B}}}
\newcommand{\AB}{\overline{\varphi}_{\mathrm{A}}}

\begin{document}
\title{Lyapunov Exponent Bounds for Random Products of Relaxed Projections}

\author{Alireza Entezari}
\address{University of Florida, Gainesville, FL 32611, USA}
\email{entezari@ufl.edu}
\author{Arunava Banerjee}
\address{University of Florida, Gainesville, FL 32611, USA}
\email{arunava@ufl.edu}
\date{}

\begin{abstract}
Current performance bounds for randomized iterated projections, such as the Kaczmarz method, are often considered tight under established per-iteration analyses, yet they are notoriously loose in practice. We derive asymptotic performance bounds that narrow this theory–practice gap by introducing a new technique for bounding the largest Lyapunov exponent of random products of relaxed projections. In addition to providing sharper iteration-complexity guarantees, the asymptotic analysis uncovers and quantifies the previously unexplained role of overrelaxation in accelerating randomized Kaczmarz, addressing a longstanding open question raised by Strohmer and Vershynin in 2007.

These results also apply directly to randomized Gauss–Seidel for symmetric positive definite systems. In contrast to deterministic SOR, where explicit formulas for optimal relaxation are available primarily under special structural assumptions, such as consistent ordering, our asymptotic analysis provides, in the randomized setting, computable estimates of the overrelaxation parameter and convergence improvements for general symmetric positive definite systems, using spectral information from the system matrix.
\end{abstract}
\maketitle

\section{Introduction}\label{sec:intro}
Alternating projections provide a fundamental framework for a broad class of iterative methods in computational mathematics~\cite{vonneumann1950functional, halperin1962product, bauschke1996projection, deutsch1997rate}, with connections to subspace correction methods underlying Gauss--Seidel, domain decomposition, and multigrid~\cite{xu2002method}. In projection methods such as Kaczmarz and Gauss--Seidel~\cite{saad2003iterative, varga1962iterative, hackbusch2016iterative}, the error is successively reduced by alternating among projection steps drawn from a set of orthogonal projectors $\{\proj_i\}_{i=1}^m$ in $\reel^n$. With relaxation parameter $\omega\in(0,2)$, for an initial error $\berr_0$, the error $\berr_k$ at the $k^{\rm th}$ step satisfies the recurrence 
\[\berr_k=(\I-\omega\proj_{i_k})\berr_{k-1},
\]
where $i_k\in\{1,\ldots,m\}$ indexes the projector selected at that step. 

Under randomization, the design and performance of these methods can be studied from a dynamical systems perspective through random products of matrices. In this work, we establish a new technique for bounding the asymptotic convergence factor of randomized projection methods, defined for each $\omega\in(0,2)$ by
\begin{equation}\label{eq:factor}
\varphi(\omega) := \lim_{k\to\infty} \norm{(\I-\omega\proj_{i_k})\cdots(\I-\omega\proj_{i_1})}^{2/k},
\end{equation}
where the indices $i_k$ are drawn independently from a fixed probability distribution $\{p_i\}_{i=1}^m$. By the Furstenberg--Kesten theorem~\cite{furstenberg1960products}, the limit\footnote{Unless otherwise specified, $\|\cdot\|$ denotes the Euclidean norm for vectors and the spectral norm for matrices. In finite dimensions, the limit is independent of the choice of matrix norm.} exists almost surely and is nonrandom.
The asymptotic convergence factor $\varphi(\omega)$ is related to the largest Lyapunov exponent of the random matrix products, which is given by $\frac12\log\varphi(\omega)$. Although Lyapunov exponents provide a fundamental characterization of the long-term behavior of random dynamical systems, they are notoriously difficult to compute or bound in general~\cite{crisanti2012products,blondel2005computationally,tsitsiklis1997lyapunov, altschuler2020lyapunov}.

The standard convergence analysis of randomized projection methods instead relies on characterizing the contraction of the {\em expected} squared error at each iteration~\cite{strohmer2009randomized, nesterov2012efficiency, wright2022optimization}. This leads to the following per-step analysis:
\begin{equation}\label{eq:err_b_bound}
\expect{\norm{(\I-\omega\proj)\berr}^2} = \norm{\berr}^2-\omega(2-\omega)\berr^T\EP\berr\leq \bigl(1-\omega(2-\omega)\mu\bigr)\norm{\berr}^2,
\end{equation}
where $\EP$ is the expected projector and $\mu$ its smallest eigenvalue, defined by
\begin{equation}\label{eq:expected_projector}
\EP := \sum_{i=1}^m p_i\proj_i,
\qquad
\mu := \lambda_{\min}(\EP).
\end{equation}
The bound $1-\omega(2-\omega)\mu$ in~\eqref{eq:err_b_bound} is minimized at $\omega=1$, suggesting that unrelaxed orthogonal projections are the optimal choice with the best guarantees on convergence in this per-step analysis. For this choice of relaxation parameter, independence of the random choices $i_k$ yields the established iteration-complexity estimate~\cite{strohmer2009randomized, leventhal2010randomized, nesterov2012efficiency} for the expected squared error,
\begin{equation*}
\expect{\norm{\berr_k}^2} \leq (1-\mu)^k\norm{\berr_0}^2,
\end{equation*}
which holds for every $k\geq 0$.

While these bounds are tight in the sense that they
are attained in certain (reducible) problems, they
are notoriously loose in practice: their per-iteration
nature often leads to conservative estimates
that substantially underestimate the observed speed
of convergence~\cite{agaskar2014randomized,
tu2017breaking, BAI2018252}. Moreover, the per-step analysis paradoxically suggests that relaxation can only hurt convergence, even though it is well known empirically that, in randomized settings as in the deterministic realm, overrelaxation can significantly improve performance. Explaining and quantifying the role of relaxation in the randomized setting has remained an open problem~\cite{strohmer2009randomized}, with roots in the deterministic setting tracing back to successive over-relaxation (SOR).

We show that asymptotic analysis bridges the longstanding gap between theoretical convergence bounds and observed performance. Despite the general difficulty of characterizing Lyapunov exponents, we show that for randomized projection methods, an effective bound on the asymptotic convergence factor $\varphi(\omega)$ can be obtained by solving a $2\times 2$ eigenvalue problem involving the {\em two smallest eigenvalues} of the expected projector $\EP$ in~\eqref{eq:expected_projector} and a single second-order statistic $\xi$ of the individual projectors, evaluated along an eigenvector associated with the smallest eigenvalue. The resulting asymptotic bound addresses both limitations of the classical per-iteration analysis: it yields substantially sharper convergence guarantees and, for the first time, quantifies the improvement in these guarantees achievable through overrelaxation.

For simplicity of notation, we assume that the smallest eigenvalue $\mu$ of the expected projector $\EP$ is positive and simple.
\begin{theorem}[Random Product of Relaxed Projections]\label{thm:RAP}
Let $\{\proj_i\}_{i=1}^m$ be a collection of orthogonal projectors, and let the indices $i_k$ be drawn independently from a fixed probability distribution. For every $\omega\in[0,2]$, the asymptotic convergence factor satisfies
\begin{equation}\label{eq:a_bound}
\varphi(\omega) = \lim_{k\to\infty}\norm{(\I-\omega\proj_{i_k})\cdots(\I-\omega\proj_{i_1})}^{2/k} \leq 1-\omega\lambda_{\min}(\mat{B}-\omega\mat{C})%
\end{equation}
almost surely, where
\begin{equation*}
\mat{B} := \begin{bmatrix}
\mu+\mu^\prime & 0 \\
0 & 2\mu
\end{bmatrix},
\qquad
\mat{C} := \begin{bmatrix}
\gamma & \sqrt{\gamma\xi} \\
\sqrt{\gamma\xi} & \xi
\end{bmatrix}.
\end{equation*}
Here, $\mu$ and $\mu^\prime$ denote the smallest and second smallest eigenvalues of $\EP$, respectively; $\xi:=\expect{(\vect{u}^T\proj\vect{u})^2}$, where $\vect{u}$ is a normalized eigenvector of $\EP$ corresponding to $\mu=\expect{\vect{u}^T\proj\vect{u}}$; and $\gamma:=\frac12(\mu+\mu^\prime)(1-\xi/\mu)$.
\end{theorem}

\begin{remark}[Rank-Agnostic Products of Projections]
The theorem is agnostic to the rank of projectors as well as the choice of probabilities with which they are chosen. Rank-1 projectors appear in randomized Kaczmarz and Gauss-Seidel/coordinate descent, while higher and mixed rank projectors appear in block methods~\cite{nesterov2012efficiency, nesterov2017efficiency, tu2017breaking} and more generally in subspace correction methods~\cite{xu2002method}.
\end{remark}

\begin{remark}[Asymptotic vs. Per-iteration bounds]
The asymptotic bound established in~\eqref{eq:a_bound} and the basic per-iteration bound from~\eqref{eq:err_b_bound}, which we refer to as the {\bf A-bound} and {\bf B-bound}, respectively, satisfy
\begin{equation}\label{eq:A_B_bound}
\varphi(\omega)
\leq \AB(\omega) := 1-\omega\lambda_{\min}(\mat{B}-\omega\mat{C})
\leq \BB(\omega) := 1-\omega(2-\omega)\mu.
\end{equation}
The two bounds coincide for certain degenerate (reducible) problems, such as when the projectors are pairwise orthogonal. More generally, as the conditioning of the problem worsens, the A-bound can improve significantly over the B-bound. The extent of this improvement depends on the gap between the eigenvalues $\mu^\prime$ and $\mu$ of the expected projector $\EP$. %
\end{remark}

\begin{figure}[htbp]
    \centering
    \includegraphics[width=\linewidth]{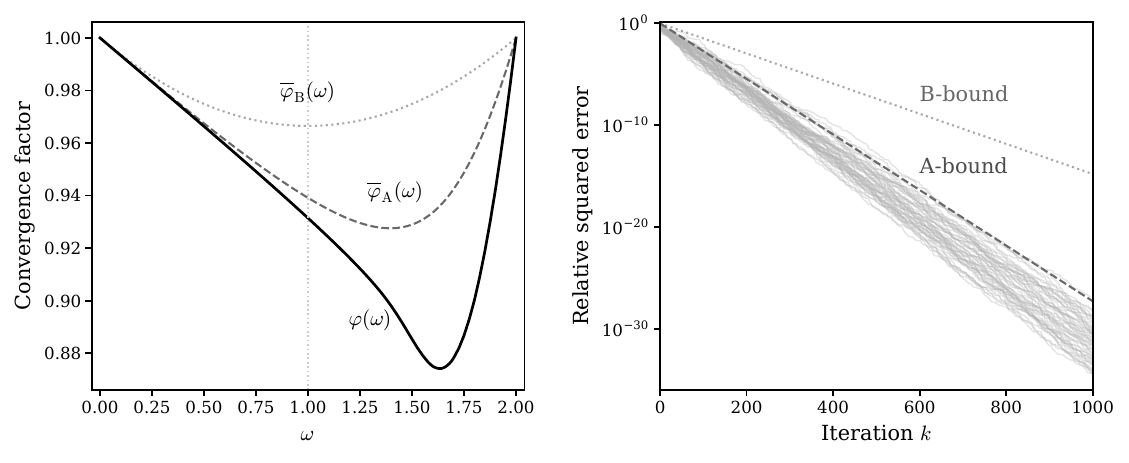}
    \caption{Convergence of eight uniformly sampled mixed-rank orthogonal projectors in $\mathbb{R}^5$. The left panel compares the asymptotic convergence factor $\varphi(\omega)$ with its A- and B-bounds. The solid curve is obtained from finite-time estimates of the top Lyapunov exponent. The right panel shows the relative squared error, on a logarithmic scale, from 50 independent trajectories at $\omega=1$, together with reference slopes corresponding to the convergence guarantees $\AB(1)$ and $\BB(1)=1-\mu$.
}
    \label{fig:mixed-rank}
\end{figure}

\begin{example}
To illustrate the bounds in~\eqref{eq:A_B_bound}, we randomly generate eight orthogonal projectors in $\mathbb{R}^5$, comprising four rank-one, three rank-two, and one rank-three projectors. Each projector is constructed by orthonormalizing the columns of an independent Gaussian matrix, and the projectors are sampled uniformly at random.

The left panel of~\FG{fig:mixed-rank} compares the asymptotic convergence factor $\varphi(\omega)$ with its A- and B-bounds. In this realization, the A-bound at $\omega=1$ improves the guaranteed contraction by a factor of approximately $1.82$ over the classical B-bound. The right panel illustrates this improvement through the relative squared-error trajectories from 50 independent runs at $\omega=1$. Their long-run slopes are substantially closer to the A-bound than to the classical B-bound.

The left panel also demonstrates the benefits of overrelaxation. The A-bound predicts an optimal relaxation parameter $\omega_\star\approx1.40$, yielding a further $19.0\%$ improvement in guaranteed contraction relative to the unrelaxed A-bound and an overall improvement by a factor of approximately $2.16$ compared with the classical B-bound. The estimated asymptotic convergence factor attains its minimum near $\omega=1.63$, illustrating that the A-bound captures the benefit of overrelaxation, although its predicted optimal parameter is conservative.
\end{example}

\begin{remark}[Consequences for error trajectories]
The bounds of~\Thm{thm:RAP} provide convergence guarantees for individual error trajectories. For any fixed nonzero initial error $\vect{z}_0$, the iterates $\vect{z}_k=(I-\omega\proj_{i_k})\vect{z}_{k-1}$ satisfy
\[
\lim_{k\to\infty}\left(\frac{\|\vect{z}_k\|^2}{\|\vect{z}_0\|^2}\right)^{1/k}\leq\AB(\omega),
\]
with probability $1$. Moreover, Proposition~\ref{prop:finite-iteration} establishes the following finite-iteration guarantee for an $n$-dimensional error space:
\[
\mathbb{E}\|\vect{z}_k\|^2\leq
\sqrt{n}\,\AB(\omega)^k\|\vect{z}_0\|^2.
\]
Thus, the asymptotic bound also yields finite-iteration guarantees that can be substantially sharper for sufficiently large $k$, despite the dimension-dependent prefactor.
\end{remark}

\begin{remark}[Sharper convergence bounds without relaxation]
Even without relaxation ($\omega=1$), the asymptotic analysis yields
\begin{equation*}
\varphi(1)\leq\AB(1)
=1-\lambda_{\min}(\mat{B}-\mat{C})
\leq 1-\mu.
\end{equation*}
This improvement narrows the theory--practice gap discussed earlier, providing sharper convergence guarantees even for unrelaxed randomized projections.
\end{remark}

\begin{remark}[Benefits of overrelaxation]
The relaxation parameter $\omega_\star \in [1,2)$ minimizing the A-bound $\AB(\omega)$ is available in closed-form and provides a quantitative guarantee for the benefits of overrelaxation.
\end{remark}

\begin{figure}[ht!]
\centering\vspace{-.25cm}
\myfig{.48}{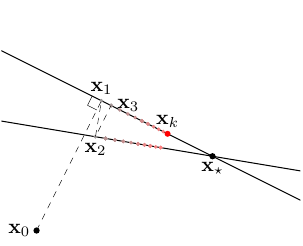}\myfig{.48}{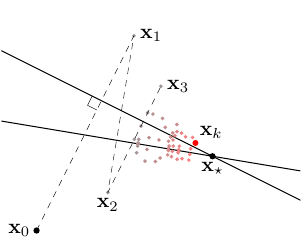}%
\caption{Alternating projections. At each iteration, one of the $m=2$ projectors in $\{\proj_1,\proj_2\}$ is selected at random. The left panel shows the unrelaxed trajectory ($\omega=1$), which minimizes the norm of the error $\berr_k:=\bx_k-\bx_\star$ at each step. The right panel illustrates how overrelaxation ($\omega=1.5$) leads to faster convergence in the long run.}
\label{fig:twolines}
\end{figure}

\begin{remark}[Computability of the bound]
The key feature of the bound is that its ingredients, which satisfy $\mu^2\leq\xi\leq\mu\leq\mu^\prime$, depend only on the distribution of individual projectors, rather than their products. Specifically, $\mu$ and $\mu^\prime$ are eigenvalues of the expected projector $\EP$, while $\xi$ is a second moment under the same sampling distribution. Thus, the asymptotic behavior of an infinite random product is bounded using only first- and second-order statistics of a single random projector. For randomized Gauss--Seidel with uniform probability distribution on a symmetric positive definite $\mA$ with unit diagonal, $\mu$ and $\mu^\prime$ are the two smallest eigenvalues of $\mA/\operatorname{tr}(\mA)$, and $\xi=n\mu^2\|\vect{u}\|_4^4$, where $\vect{u}$ is a normalized eigenvector of $\mA$ corresponding to $\mu$. For randomized Kaczmarz with uniform sampling on a matrix $\mA$ with unit row norms, $\mu$ and $\mu^\prime$ are the two smallest squared singular values of $\mA/\|\mA\|_F$, and $\xi=m\mu^2\|\vect{v}\|_4^4$, where $\vect{v}$ is the unit left singular vector  of $\mA$ corresponding to its smallest singular value. Thus, the ingredients are accessible from standard spectral quantities and, in the context of statistical learning~\cite{defossez2015averaged} or streaming applications, can also be estimated online.
\end{remark}

\begin{remark}[Degenerate eigenvalues]
The assumption that the smallest eigenvalue $\mu$ of
$\EP$ is positive and simple is made for convenience
of notation. When $\EP$ is singular, the iterations
converge to the common nullspace of the
projectors. The asymptotic convergence factor
toward this subspace can be characterized by
the same bound, applied on its orthogonal
complement, taking $\mu$ to be the smallest
positive eigenvalue of $\EP$. Moreover, repeated eigenvalues are nongeneric. When $\mu$ has multiplicity greater than one,
an analogous bound can be obtained by considering
the entire eigenspace associated with $\mu$,
with an appropriate second-order statistic,
as discussed in~\SC{sec:discussion}.
\end{remark}

\begin{example}[Practical Significance]
\label{ex:iteration-complexity}
The sharper asymptotic bounds translate into guarantees of convergence to a prescribed tolerance within substantially fewer iterations. As an illustration, we apply randomized Gauss–Seidel to a Gaussian-generated symmetric positive definite system of dimension $n=250$. For a relative expected squared-error tolerance of $10^{-8}$, the classical per-iteration B-bound guarantees convergence within approximately $2.10\times10^9$ coordinate updates. For the same algorithm, sampling distribution, and relaxation parameter $\omega=1$, the finite-iteration estimate derived from the A-bound reduces this guarantee to $1.21\times10^9$ updates, a reduction of $42.5\%$ achieved solely through sharper analysis. Moreover, the asymptotic analysis identifies further improvements through overrelaxation, which are not captured by the classical per-iteration bound. Choosing the predicted optimal relaxation parameter $\omega\approx1.62$ in this realization reduces the guaranteed iteration count to $7.47\times10^8$, representing $64.5\%$ fewer iterations than the classical guarantee, or a reduction by a factor of approximately $2.82$. Thus, the new analysis provides substantially improved iteration-complexity guarantees, both for the original algorithm and through overrelaxation.
\end{example}

Further examples illustrating applications to randomized Kaczmarz, gossip algorithms, and block subspace corrections are presented in~\SC{sec:app}.

The main technical contribution of this work is a new approach for bounding the spectral radius of a class of positive linear maps (superoperators) arising from random products of relaxed orthogonal projections. By exploiting the spectral structure of these operators, we reduce a high-dimensional eigenvalue problem to a $2\times2$ eigenvalue problem involving only the two smallest eigenvalues of the expected projector and a single second-moment quantity. Combined with the connection between second-moment growth and the largest Lyapunov exponent, this yields the asymptotic convergence bound in Theorem~\ref{thm:RAP}.

The remainder of the paper is organized as follows. In \SC{sec:framework}, we develop the second-moment superoperator framework and establish its connection to the asymptotic convergence factor. In \SC{sec:radius}, we introduce a new technique for bounding the spectral radius of the second-moment superoperator and prove Theorem~\ref{thm:RAP}. In \SC{sec:app}, we apply the resulting bounds to randomized Kaczmarz, Gauss--Seidel, and block subspace correction methods, quantifying improvements over classical per-iteration guarantees and the acceleration achievable through overrelaxation. Finally, \SC{sec:discussion} examines related work, implications, and possible extensions.

\section{Analytical Framework}\label{sec:framework}

\paragraph{Notation}
We use bold lowercase characters to denote vectors (e.g., $\bx \in \reel^n$ and $\bb \in \reel^m$), and bold uppercase characters to denote matrices (e.g., $\mA \in \reel^{m\times n}$). The eigenvalues of a symmetric matrix $\mat{M}$ are indexed in increasing order:
\begin{equation*}
    \lambda_{\min}(\mat{M}) := \lambda_1(\mat{M})
    \le \lambda_2(\mat{M}) \le \cdots
    \le \lambda_{\max}(\mat{M}).
\end{equation*}
The corresponding eigenvectors are denoted by $\vect{v}_i(\mat{M})$.
Calligraphic letters denote linear operators acting on the space of $n\times n$ matrices (superoperators). We denote their eigenvalues by $\lambda_i(\T)$ and corresponding eigenvectors (eigenmatrices) by $\mat{V}_i(\T)$.
The Frobenius inner product is denoted by
$\langle \mat{M}, \mat{M}^\prime \rangle := \tr(\mat{M}^T \mat{M}^\prime)$,
with associated norm $\|\cdot\|_F$.
For a superoperator $\T$, its spectral radius is
$\rho(\T):=\max_i |\lambda_i(\T)|$.

\subsection{Second-Moment Superoperator}
As discussed in~\SC{sec:intro}, in randomized alternating projections the error evolves as a recurrence
$\berr_{k}=(\I-\omega\proj_{i_{k}})\berr_{k-1}$,
where the projectors are selected independently by choosing $i_k$ from a fixed probability distribution. One can analyze the evolution of the second-moment matrix $\mat{\Sigma}_k:=\expect{\berr_k\berr_k^T}$, starting from $\mat{\Sigma}_0:=\berr_0\berr_0^T$. For $k>0$, independence of the random choices $i_k$ implies that the second-moment matrix evolves according to the recurrence
\[
\mat{\Sigma}_{k}=\expect{(\I-\omega\proj)\mat{\Sigma}_{k-1}(\I-\omega\proj)^T},
\]
where the expectation is taken over random choice from projectors $\{\proj_i\}_{i=1}^m$, with $i\in\{1,\ldots,m\}$ drawn from the fixed probability distribution $\{p_i\}_{i=1}^m$ introduced in~\SC{sec:intro}. This evolution is governed by a linear map (superoperator) acting on the space of $n \times n$ matrices
\[
    \A(\mat{X})
    := \expect{(\I-\omega\proj)\mat{X}(\I-\omega\proj)^T}.
\]
Consequently, the second-moment matrix evolves through repeated applications of the superoperator, $\mat{\Sigma}_k=\A^k(\mat{\Sigma}_0)$. The superoperator admits a Kronecker-product representation, and its spectral radius, which we show in the next subsection provides an upper bound on the asymptotic convergence factor, is denoted by $\overline{\varphi}(\omega)$:
\begin{equation}\label{eq:covariance_dynamics}
    \A := \expect{(\I-\omega\proj)\otimes(\I-\omega\proj)},
    \qquad
    \overline{\varphi}(\omega) := \rho(\A).
\end{equation}

\subsection*{Perron-Frobenius Theory} Analyzing power iterations on $\mat{\Sigma}_{k+1} = \A(\mat{\Sigma}_k)$, suggests that if $\mat{\Sigma}_k$ were to converge in $\reel^{n \times n}$ to $\mat{0}$ along an eigenvector of the  superoperator $\A$, that eigenvector needs to be positive (semi) definite since $\mat{\Sigma}_k \succcurlyeq \mat{0}$. In fact, the generalization of the Perron-Frobenius theorem on $C^*$-algebras~\cite{evans1978spectral, farenick1996irreducible}, that holds for $\A$ (see Appendix~\ref{sec:PF}), shows that the spectral radius of $\A$ is attained by an eigenvalue whose corresponding eigenvector can be chosen to be a positive semi-definite matrix that describes the asymptotic evolution of the $\mat{\Sigma}_k$. Therefore, the spectral radius $\overline{\varphi}(\omega)$ is an eigenvalue of $\A$. Furthermore, if the set of projectors $\{\proj_i\}_{i=1}^m$ is {\em irreducible} (i.e., has no nontrivial common invariant subspace), then the theory guarantees that (i) the eigenvector of $\A$ corresponding to the spectral radius is positive definite, (ii) it is the only eigenvector that is non-negative, and (iii) that the eigenvalue that attains the spectral radius is simple (not repeated). This implies $\overline{\varphi}(\omega) = \lambda_{\max}(\A)$ and that it is a smooth function of $\omega \in (0, 2)$.

\subsection{Connection with the Asymptotic Convergence Factor}
The spectral radius of the second-moment superoperator $\A$ provides a classical Jensen-type bound on the Lyapunov exponent and, consequently, on the asymptotic convergence factor $\varphi(\omega)$ defined in~\eqref{eq:factor}. Such bounds are commonly associated with annealed analysis~\cite{crisanti2012products,agaskar2014randomized}.

\begin{proposition}\label{prop:lyapunov_jensen}
The asymptotic convergence factor is bounded by
$\varphi(\omega)\leq\overline{\varphi}(\omega)$.
\end{proposition}

\begin{proof}
For a nonzero initial error $\berr_0$ in general position, let
$\mat{\Sigma}_0=\berr_0\berr_0^T$.
Since $\mat{\Sigma}_k=\A^k(\mat{\Sigma}_0)$, using operator norm we get $\norm{\A^k}  \ge \norm{\mat{\Sigma}_k} /  \norm{\mat{\Sigma}_0}$. 
Taking the $k^\text{th}$ root, then limit superior on both sides and applying Gelfand's formula then gives
\[
\rho(\A)
=\lim_{k\to\infty}\norm{\A^k}^{1/k}
\geq \limsup_{k\to\infty}\norm{\mat{\Sigma}_k}^{1/k},
\]
since the fixed factor $\norm{\mat{\Sigma}_0}^{1/k}$
limits to one.

While Gelfand's formula holds for any norm on $\mat{\Sigma}_k$, since it is positive semidefinite, its trace
norm provides $\tr(\mat{\Sigma}_k)=\expect{\norm{\berr_k}^2}$.
Consequently, Jensen's inequality yields
\[
\overline{\varphi}(\omega)=
\rho(\A)
\geq \limsup_{k\to\infty}
\left(\expect{\norm{\berr_k}^2}\right)^{1/k}
\geq \limsup_{k\to\infty}
\expect{\norm{\berr_k}^{2/k}}.
\]

By Furstenberg-Kesten~\cite{furstenberg1960products} and the more general multiplicative ergodic theorem of Oseledets~\cite{oseledec1968multiplicative},
the asymptotic convergence factor is nonrandom and attained
almost surely for an initial error in general position
(see also~\cite[Eq.~(2.8)]{crisanti2012products}):
\[
\lim_{k\to\infty}\norm{\berr_k}^{2/k}
=\varphi(\omega)
\qquad\text{almost surely}.
\]
Applying Fatou's lemma to the nonnegative random variables
$\norm{\berr_k}^{2/k}$ gives
\[
\begin{aligned}
\varphi(\omega) = \expect{\liminf_{k\to\infty}\norm{\berr_k}^{2/k}}
\leq \liminf_{k\to\infty}
\expect{\norm{\berr_k}^{2/k}}
\leq \limsup_{k\to\infty}
\expect{\norm{\berr_k}^{2/k}}
\leq \overline{\varphi}(\omega).
\end{aligned}
\]
\end{proof}

While the analysis of second-moment dynamics has a long history in stochastic approximation~\cite{solo1989limiting, aguech2000perturbation}, with the superoperator perspective developed in~\cite{murata1999statistical} and subsequently applied to error analysis in~\cite{agaskar2014randomized, defossez2015averaged}, obtaining performance bounds sharper than those derived from per-iteration analysis has remained challenging. A central difficulty is relating the spectral radius of the superoperator $\A$ to the spectral properties of the original problem (e.g., the spectrum of $\mA$ in Kaczmarz or Gauss--Seidel). The A-bound $\AB(\omega)$ introduced in~\Thm{thm:RAP} addresses this challenge by bounding the spectral radius of $\A$ using only statistics of individual random projectors, without analyzing products of random projectors. To establish~\Thm{thm:RAP}, we develop a new technique for bounding the spectral radius of $\A$ for every $\omega\in[0,2]$ and prove that $\overline{\varphi}(\omega)\le\AB(\omega)$.

\section{Bounding the Spectral Radius of the Second-Moment Superoperator}
\label{sec:radius}

The Kronecker representation of the second-moment superoperator
$\A$ in~\eqref{eq:covariance_dynamics} yields the decomposition
\begin{equation}\label{eq:BC}
    \A = \IxI - \omega \B + \omega^2 \C
       = \IxI - \omega(\B - \omega\C),
\end{equation}
where $\IxI := \I \otimes \I$,
$\B := \EP \oplus \EP = \I \otimes \EP + \EP \otimes \I$,
and $\C := \EPxP$.
The expectations are taken over the collection of orthogonal
projectors $\{\proj_i\}_{i=1}^m$, selected independently with
probabilities $\{p_i\}_{i=1}^m$, as introduced in
~\SC{sec:intro}. The superoperator $\B$ depends only
on the first-order statistics of the random projector,
through $\EP$, whereas $\C$ captures second-order statistics
through $\EPxP$. Both are independent of the relaxation
parameter $\omega$.

While $\B$ is completely determined by the expected projector
$\EP$, the superoperator $\C$ contains additional information
about the distribution of individual projectors. Moreover,
$\B$ and $\C$ generally do not share an eigenbasis, so the
spectral radius of $\A$ depends on their interaction rather
than their individual spectra alone. Our objective is to
bound the spectral radius of $\A$ using only limited information
about $\C$. We show that, together with $\EP$, the single
second-order statistic $\xi$ introduced in
Theorem~\ref{thm:RAP} suffices to construct a rank-one
surrogate superoperator $\C^\star$ that provides an effective
upper bound on the spectral radius of $\A$ for every
$\omega\in[0,2]$.

By the Perron--Frobenius theory for positive linear maps,
the spectral radius of $\A$ is attained by an eigenvalue.
This implies $\rho(\A)=\lambda_{\max}(\A)$. Consequently,
the decomposition in~\eqref{eq:BC} gives
\begin{equation}\label{eq:AtoBC}
    \rho(\A)=\lambda_{\max}(\A)
    =1-\omega\lambda_{\min}(\B-\omega\C).
\end{equation}
Thus, obtaining an upper bound on the spectral radius of
$\A$ amounts to establishing a lower bound on the smallest
eigenvalue of $\B-\omega\C$, which we interpret geometrically
as the minimum separation between the quadratic forms of
$\B$ and $\omega\C$.

Classical perturbation theory provides tools for studying
how the eigenvalues of $\B-\omega\C$ deviate from those
of $\B$ as $\omega$ varies. In particular, Hadamard
variation formulae can be used to bound these deviations
through eigenvalue derivatives~\cite{tao2011universality}.
However, such general perturbation bounds offer only
marginal improvements over the classical B-bound.
Instead, we exploit the structural relationship between
$\B$ and $\C$ to develop a geometric approach for
bounding the smallest eigenvalue of $\B-\omega\C$.

A fundamental structural property, established in
Lemma~\ref{lem:loewner}, is the Loewner inequality
$\B\succeq 2\C$, with $(\B-2\C)(\I)=0$.
Consequently, for every $\omega\in[0,2]$, the quadratic
form of $\omega\C$ lies below that of $\B$, and the two
coincide along the identity direction at $\omega=2$.

To develop intuition behind our approach, we adopt a
geometric view of Rayleigh quotients over the unit sphere
in the space of matrices equipped with the Frobenius
inner product. For any self-adjoint superoperator
$\T$, we denote its Rayleigh quotient by
$R_{\T}(\mat{V})
:=\langle\mat{V},\T(\mat{V})\rangle/
\|\mat{V}\|_F^2$ for $\mat{V}\ne\mat{0}$.
Let $\mathbb{S}:=\{\mat{V}\in\reel^{n\times n}:
\|\mat{V}\|_F=1\}$ denote the unit sphere.
On $\mathbb{S}$, the Rayleigh quotient agrees with
the corresponding quadratic form. By the variational
characterization of the smallest eigenvalue,
\begin{equation}\label{eq:variational}
    \lambda_{\min}(\B-\omega\C)
    =\min_{\mat{V}\in\mathbb{S}}
    \left\{
        R_{\B}(\mat{V})-\omega R_{\C}(\mat{V})
    \right\}.
\end{equation}

For the illustration in~\FG{fig:B_growing_C},
we consider projectors in $\reel^2$ and restrict
$\mat{V}$ to symmetric $2\times2$ matrices, which form
a three-dimensional vector space. For each unit matrix
$\mat{V}$, the values $R_{\B}(\mat{V})$ and
$\omega R_{\C}(\mat{V})$ specify the scaling of
$\mat{V}$, producing the black and red surfaces,
respectively. The smallest eigenvalue corresponds
to the minimum separation between the two surfaces
measured along the same direction on the unit sphere.
At $\omega=0$, the minimizing direction is $\mU$,
with minimum separation $2\mu$. As $\omega$ increases,
the minimum separation decreases, reaching zero at
$\omega=2$, where the normalized identity $\I/\sqrt{n}$
is a minimizing direction.

\begin{figure}[t]
    \centering
    \includegraphics[width=0.75\linewidth]
        {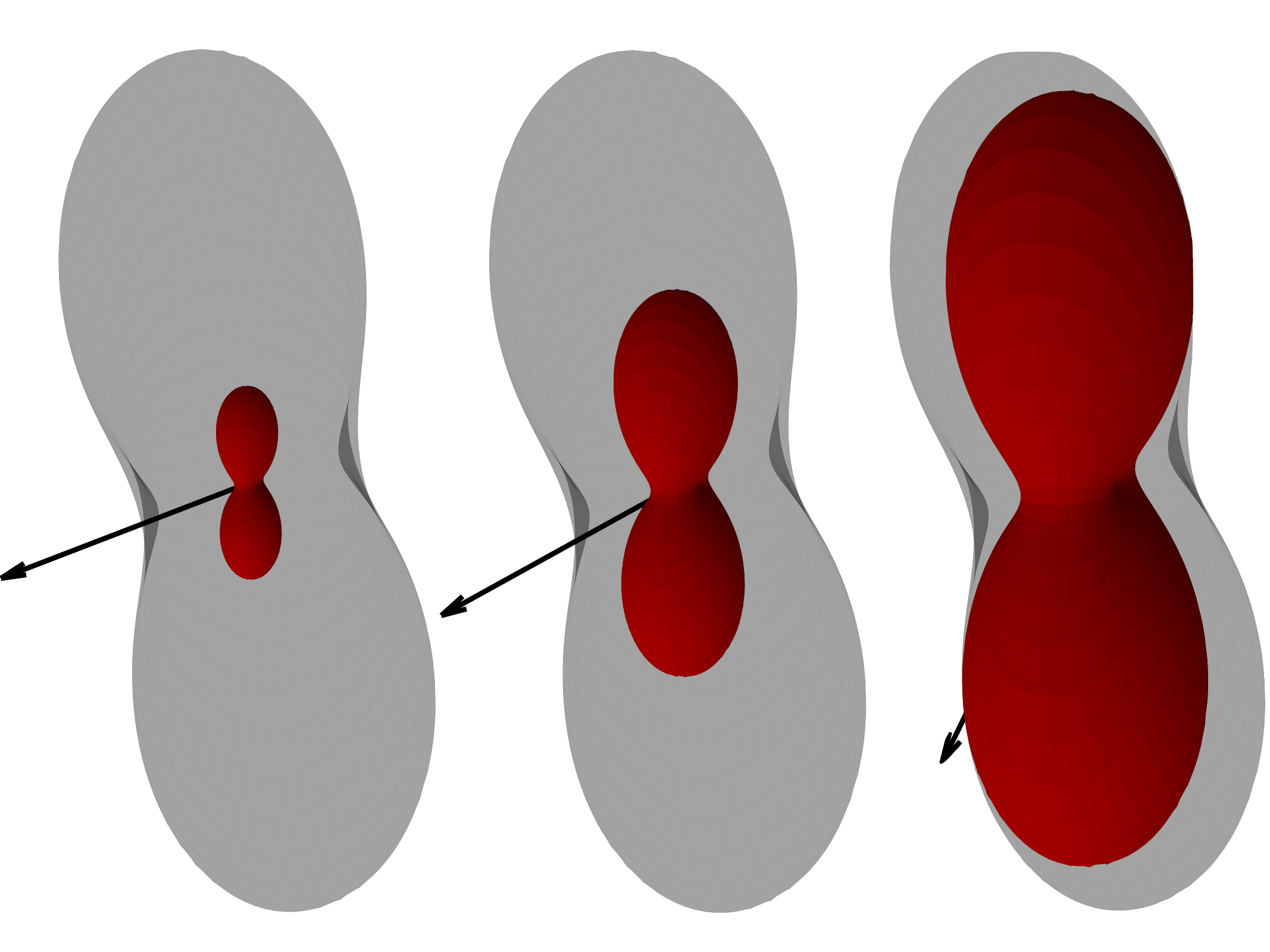}
    \caption{Radial plots of the Rayleigh quotients
    associated with superoperators $\B$ (black) and
    $\omega\C$ (red) over the unit sphere. By the
    variational characterization, the smallest eigenvalue
    $\lambda_{\min}(\B-\omega\C)$ corresponds to the
    minimum separation between the two surfaces along
    a common direction, identified by the eigenvector
    $\mat{V}_1(\B-\omega\C)$. The minimizing direction
    is shown for $\omega=1/2$ (left), $\omega=1$
    (middle), and $\omega=2$ (right).}
    \label{fig:B_growing_C}
\end{figure}

While the quadratic forms are defined over
$\reel^{n\times n}$, their geometric properties
can be investigated through lower-dimensional slices.
For a nonzero subspace $S$, let $\mathcal{P}_S$
denote the orthogonal projector onto $S$ and define
the {\em compression} of a superoperator $\T$
to $S$ by $\T_S:=\mathcal{P}_S\T|_S:S\rightarrow S$.
The Rayleigh quotients of $\T$ and
$\T_S$ agree on $S$, namely
$R_{\T_S}(\mat{V})=R_{\T}(\mat{V})$
for every nonzero $\mat{V}\in S$. Consequently,
\begin{equation}\label{eq:variational_slice}
    \lambda_{\min}(\B_S-\omega\C_S)
    =\min_{\mat{V}\in S\cap\mathbb{S}}
    \left\{
        R_{\B}(\mat{V})-\omega R_{\C}(\mat{V})
    \right\} \geq\lambda_{\min}(\B-\omega\C).
\end{equation}
In an orthonormal basis of a two-dimensional subspace,
these compressed quadratic forms are represented by
$2\times2$ symmetric matrices, and the smallest
eigenvalue of their difference gives the minimum
separation within that slice. Importantly, equality
holds in~\eqref{eq:variational_slice} whenever $S$
contains a minimizing eigenvector of $\B-\omega\C$.

\subsection*{Proof strategy.}
We exploit the spectral structure of $\B$ to identify
a distinguished two-dimensional subspace
$S^\star=\operatorname{span}\{\mU^\prime,\mU\}$,
where $\mU$ and $\mU^\prime$ are eigenvectors corresponding to its
smallest and second smallest eigenvalues, respectively. Using this subspace, we construct
a rank-one surrogate $\C^\star$ and establish
$\lambda_{\min}(\B-\omega\C^\star)
\leq\lambda_{\min}(\B-\omega\C)$ for every $\omega\in[0,2]$, {\em without} requiring $\C^\star\succeq\C$ in the Loewner order (see \FG{fig:eclipse_c_star}).

Our proof proceeds through comparisons of quadratic
forms on two-dimensional subspaces containing $\mU$.
For any unit matrix $\mat{V}$ orthogonal to $\mU$,
we consider the subspace
$S_{\mat{V}}=\operatorname{span}\{\mat{V},\mU\}$
and construct a corresponding rank-one surrogate
$\C^\star_{\mat{V}}$.
We first show that this surrogate yields a minimum
Rayleigh quotient no greater than that of the
compression of $\B-\omega\C$ to $S_{\mat{V}}$,
uniformly for all $\omega\in[0,2]$.
Comparing these surrogates across subspaces then
identifies $S^\star$ as the optimal choice,
again uniformly in $\omega$.
Finally, for each fixed $\omega\in(0,2]$, we consider
the subspace spanned by $\mU$ and a minimizing eigenvector
of $\B-\omega\C$. The minimum Rayleigh quotient on this
subspace equals $\lambda_{\min}(\B-\omega\C)$.
The two-dimensional comparisons, which hold for every
subspace and every $\omega$, therefore imply
$\lambda_{\min}(\B-\omega\C^\star)
\leq\lambda_{\min}(\B-\omega\C)$.

\subsection{Structural Properties}

We begin by examining the spectral structure of $\B$ using standard properties of Kronecker products and sums. For symmetric matrices $\mat{M}$ and $\mat{M}^\prime$, the eigenvalues of $\mat{M}\otimes\mat{M}^\prime$ are the products $\lambda_i(\mat{M})\lambda_j(\mat{M}^\prime)$, whereas those of $\mat{M}\otimes\I+\I\otimes\mat{M}^\prime$ are the sums $\lambda_i(\mat{M})+\lambda_j(\mat{M}^\prime)$. In both cases, the corresponding eigenvectors are Kronecker products of the eigenvectors of the constituent matrices.

\begin{fact}\label{prop:BC}
Applying the properties of Kronecker sums to $\B=\EP\oplus\EP$ shows that its eigenvalues are pairwise sums of the eigenvalues of $\EP$. In particular, its smallest eigenvalue is $\lambda_1(\B)=2\mu$, with corresponding normalized eigenvector $\mU:=\vect{u}\vect{u}^T$, where $\vect{u}$ is a normalized eigenvector of $\EP$ corresponding to $\mu$. Its second smallest eigenvalue is $\lambda_2(\B)=\mu+\mu^\prime$, with corresponding normalized eigenvector $\mU^\prime:=\vect{u}{\vect{u}^\prime}^T$, where $\vect{u}^\prime$ is a normalized eigenvector of $\EP$ corresponding to $\mu^\prime$. This eigenvalue has multiplicity at least two, and its eigenspace contains ${\rm span}\{\vect{u}{\vect{u}^\prime}^T,\vect{u}^\prime\vect{u}^T\}$. Moreover, $\B\succ\mat{0}$ if and only if $\mu>0$. Since $\C$ is a convex combination of Kronecker products of orthogonal projectors, $\C\succcurlyeq\mat{0}$. The superoperator $\A$ is symmetric but need not be positive semidefinite.
\end{fact}

We now establish the fundamental Loewner inequality invoked in the geometric interpretation above, which guarantees that the quadratic form of $\omega\C$ lies below that of $\B$ for $\omega\in[0,2]$, with the two surfaces touching at $\omega=2$.

\begin{lemma}[Loewner ordering of superoperators: $\B \succcurlyeq 2\C$\label{lem:loewner}]
\end{lemma}
\begin{proof}
    Let $\vect{v}_1, \dots, \vect{v}_n$ denote an eigenbasis for a projector $\proj$ with corresponding eigenvalues $\lambda_1,\dots,\lambda_n$. Then $\vect{v}_i \otimes \vect{v}_j$ is an eigenvector for $\proj \otimes \I + \I \otimes \proj$ and $\proj \otimes \proj$: 
    \begin{align*}
        \left(\proj \otimes \I + \I \otimes \proj - 2 \proj \otimes \proj\right)(\vect{v}_i \otimes \vect{v}_j) = (\lambda_i + \lambda_j - 2 \lambda_i \lambda_j) (\vect{v}_i \otimes \vect{v}_j).
    \end{align*}
    Since $\proj$ is a projection $\lambda_i \in \{0, 1\}$ then $\lambda_i + \lambda_j - 2 \lambda_i \lambda_j \in \{0, 1\}$. When $\proj$ is an orthogonal projection (symmetric) with an orthonormal eigenbasis $\vect{v}_i$, the set $\vect{v}_i \otimes \vect{v}_j$ forms an orthonormal eigenbasis for these superoperators and $\proj \otimes \I + \I \otimes \proj - 2 \proj \otimes \proj \succcurlyeq \mat{0}$. This establishes the result since expectation preserves the Loewner ordering
    \begin{align*}
        \B - 2\C = \expect{\proj \otimes \I + \I \otimes \proj - 2 \proj \otimes \proj} \succcurlyeq \mat{0}.
    \end{align*}    
    Moreover, we have $\B(\I) = 2\EP$ and $\C(\I) = \EP$. This means $\B - 2\C$ is positive semi-definite with $\lambda_1(\B - 2\C) = 0$ and a corresponding eigenvector is $\mat{V}_1(\B - 2\C) = \frac{1}{\sqrt{n}}\I$.
\end{proof}

The Loewner inequality in Lemma~\ref{lem:loewner} also reveals an important asymmetry in the effect of relaxation. As the following lemma shows, overrelaxation yields a spectral radius no larger than the corresponding underrelaxation, so a minimizer of $\overline{\varphi}(\omega)$ can be sought in $[1,2)$. To make the dependence on the relaxation parameter explicit, we write $\A(\omega):=\IxI-\omega(\B-\omega\C)$.

\begin{lemma}[Overrelaxation is no worse than the corresponding underrelaxation] \label{lem:overrelaxation}
For $t \in [0,1]$, 
$
\overline{\varphi}(1 - t) \ge \overline{\varphi}(1 + t)
$
\end{lemma}
\begin{proof}
We have $\A(1-t) \succcurlyeq \A(1+t)$ since $\A(1-t) - \A(1+t) = 2t \B - 4t \C = 2t (\B - 2\C) \succcurlyeq \mat{0}$. The Loewner order implies the order on the corresponding eigenvalues and we recall that spectral radius is attained by an eigenvalue due to Perron-Frobenius.
\end{proof}

\subsection{Spectral Bounds and Their Dependence on Relaxation}
We now examine how the smallest eigenvalue of $\B-\omega\C$, and consequently the spectral radius of $\A(\omega)$, vary with the relaxation parameter $\omega$. We first recover the classical B-bound using the Loewner inequality established in Lemma~\ref{lem:loewner}. We then investigate the concavity of the spectral gap and its endpoint derivatives, revealing the role of the second-order statistic $\xi$ in the behavior of these bounds.

\begin{lemma}[The B-bound] {$\lambda_i(\B - \omega\C) \ge (1-\omega/2)\lambda_i(\B)$} 
\label{lem:bbound}
\end{lemma}
\begin{proof}
    Since $\B \succcurlyeq 2\C$, we have $\B - \omega \C \succcurlyeq \B - \omega/2 \B = (1 - \omega/2) \B$. The Loewner order implies the order on the corresponding eigenvalues.
\end{proof}
The above lemma re-establishes the B-bound since $\lambda_{\max}(\A(\omega)) = 1 - \omega \lambda_1(\B - \omega \C) \le 1 - \omega (1 - \omega/2) \lambda_1(\B) = 1 - \omega (1 - \omega/2) 2 \mu$. This means the B-bound is a bound on the spectral radius: $\overline{\varphi}(\omega) \le \BB(\omega)$.

\begin{lemma}[Concavity and Convexity]\label{lem:concavity}
The smallest eigenvalue $\lambda_{\min}(\B-\omega\C)$ is a concave function of $\omega$, whereas $\lambda_{\max}(\A(\omega))$ is a convex function of $\omega$.
\end{lemma}
\begin{proof}
By the min--max theorem, $\lambda_{\min}(\B-\omega\C)$ is the pointwise minimum of affine functions of $\omega$, and is therefore concave. Similarly, $\lambda_{\max}(\A(\omega))$ is the pointwise maximum of the quadratic functions $\langle\mat{V},\A(\omega)(\mat{V})\rangle$ over all unit matrices $\mat{V}$. Since $\C\succcurlyeq\mat{0}$, each of these quadratic functions is convex in $\omega$, and hence so is their maximum.
\end{proof}

\begin{figure}[ht!]
\centering
\myfig{.5}{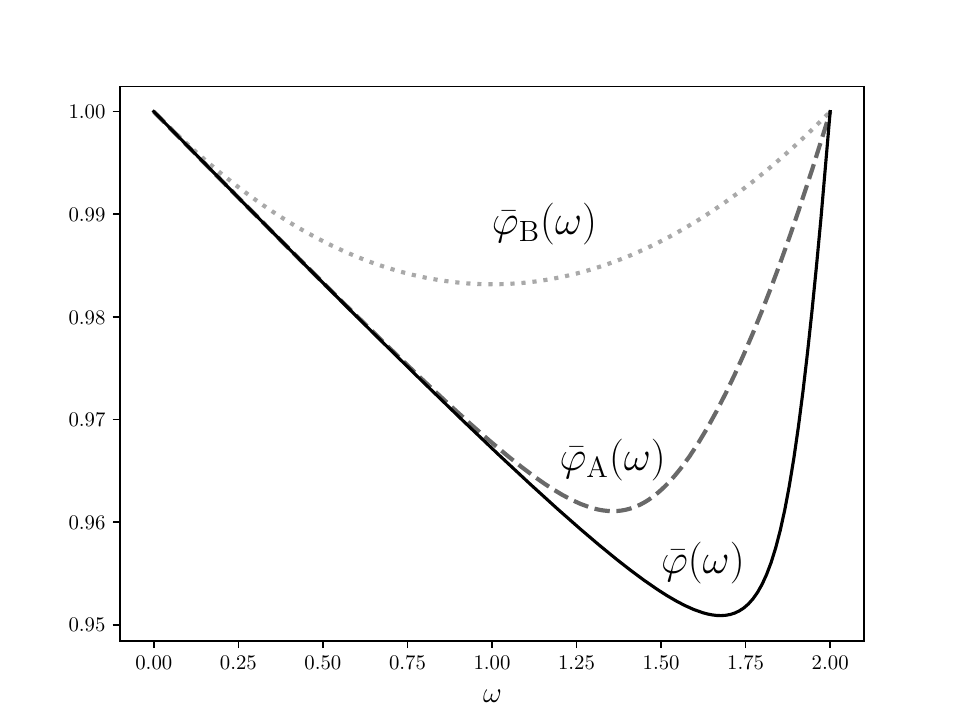}\myfig{.5}{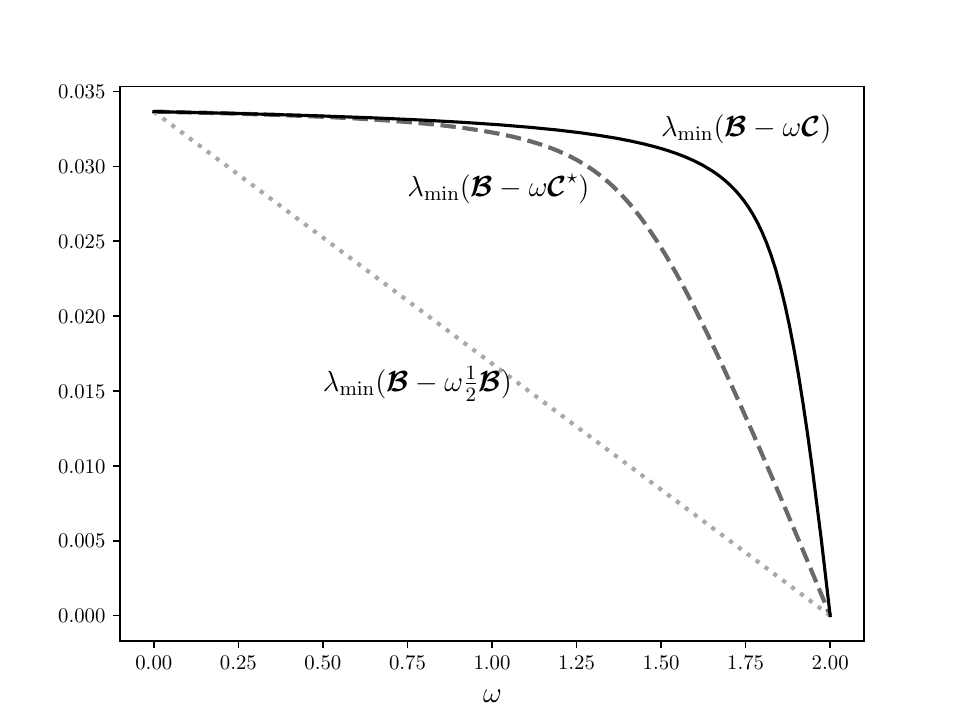}
\caption{Upper bounds on the spectral radius of $\A$ (left), for a $5\times 5$ matrix $\mA$ with random entries, correspond to lower bounds on the spectral gap of $\B - \omega \C$ (right).}
\label{fig:toy_A_vs_BC}
\end{figure}

Together, Lemmas~\ref{lem:concavity} and~\ref{lem:loewner} establish the convexity of $\overline{\varphi}(\omega)$, while the overrelaxation inequality in Lemma~\ref{lem:overrelaxation} ensures that a minimizer can be sought in $[1,2)$. Figure~\ref{fig:toy_A_vs_BC} illustrates the resulting hockey-stick shape of the spectral radius (left) and the corresponding concavity of the spectral gap (right), together with their respective bounds.

Irreducibility of the projectors $\{\proj_i\}_{i=1}^m$ guarantees that spectral radius of $\A(\omega)$ is a simple eigenvalue due to the Perron-Frobenius theorem for positive linear maps (see Appendix~\ref{sec:PF}). This means $\lambda_{\min}(\B - \omega \C)$ is also a simple eigenvalue and hence differentiable with respect to $\omega$. Since we will be focusing on the smallest eigenvalue and the corresponding normalized eigenvector of $\B - \omega \C$, for convenience of notation, we define each as a function of $\omega$:
\begin{equation}\label{eq:lam_omega}
    \lambda(\omega) := \lambda_1(\B - \omega \C), \quad {\rm and} \quad \mat{V}(\omega) := \mat{V}_1(\B - \omega \C).
\end{equation}
In particular, we have $\lambda(0) = \lambda_1(\B) = 2\mu$ and $\lambda(2) = \lambda_1(\B - 2\C) = 0$ according to Lemma~\ref{lem:loewner}.

Here is where $\xi$, the last ingredient in~\Thm{thm:RAP}, comes into play:

\begin{lemma}[Derivatives with respect to $\omega$]\label{lem:derivatives}
Suppose the smallest eigenvalue $\mu$ of $\EP$ is simple. Then the right derivative at $\omega=0$ satisfies
\begin{equation}
    \dot{\lambda}(0^+)=-\langle\mU,\C(\mU)\rangle=-\xi.
\end{equation}
If zero is a simple eigenvalue of $\B-2\C$, then the left derivative at $\omega=2$ satisfies
\begin{equation}
    \dot{\lambda}(2^-)=-\frac{1}{n}\tr{\EP}.
\end{equation}
\end{lemma}
\begin{proof}
At any $\omega$ where $\lambda(\omega)$ is simple, Hadamard's first variation formula, together with the definition of $\C$ in~\eqref{eq:BC}, gives
\begin{equation*}
    \dot{\lambda}(\omega)=-\langle\mat{V}(\omega),\C(\mat{V}(\omega))\rangle=-\expect{\tr{\left(\mat{V}(\omega)^T\proj\mat{V}(\omega)\proj\right)}}.
\end{equation*}
Since $\mu$ is simple, the smallest eigenvalue $2\mu$ of $\B$ is simple, with normalized eigenvector $\mat{V}(0)=\mU=\vect{u}\vect{u}^T$. Therefore,
\begin{equation*}
    \dot{\lambda}(0^+)=-\expect{\tr{\left(\vect{u}\vect{u}^T\proj\vect{u}\vect{u}^T\proj\right)}}=-\expect{\left(\vect{u}^T\proj\vect{u}\right)^2}=-\xi,
\end{equation*}
as defined in~\Thm{thm:RAP}. Moreover, Lemma~\ref{lem:loewner} establishes that $\B-2\C$ has a zero eigenvalue with normalized eigenvector $\I/\sqrt{n}$. Under the additional simplicity assumption, Hadamard's formula yields
\begin{equation*}
    \dot{\lambda}(2^-)=-\frac1n\left\langle\I,\C\left(\I\right)\right\rangle=-\frac{1}{n}\expect{\tr{\proj^2}}=-\frac{1}{n}\tr{\EP},
\end{equation*}
where the last equality follows from $\proj^2=\proj$.
\end{proof}

Recalling from~\eqref{eq:AtoBC} that $\overline{\varphi}(\omega)=\lambda_{\max}(\A(\omega))$, we obtain $\dot{\overline{\varphi}}(0^+)=-2\mu$ and $\ddot{\overline{\varphi}}(0^+)=2\xi$. If zero is a simple eigenvalue of $\B-2\C$, we also have $\dot{\overline{\varphi}}(2^-)=\frac{2}{n}\tr{\EP}$. By convexity and concavity, the tangent lines at the endpoints provide lower bounds for $\overline{\varphi}(\omega)$ and upper bounds for $\lambda(\omega)$, respectively, as illustrated in~\FG{fig:toy_A_vs_BC}.

When the endpoint eigenvalues are repeated, the corresponding one-sided derivatives still exist, but are determined by extremal Rayleigh quotients over their eigenspaces. Writing $E_0:=\ker(\B-2\mu\IxI)$ and $E_2:=\ker(\B-2\C)$, we have
\begin{align*}
    \dot{\lambda}(0^+)=-\max_{\substack{\mat{V}\in E_0\cap\mathbb{S}}}\langle\mat{V},\C(\mat{V})\rangle, \qquad \dot{\lambda}(2^-)=-\min_{\substack{\mat{V}\in E_2\cap\mathbb{S}}}\langle\mat{V},\C(\mat{V})\rangle.
\end{align*}
The formulas in Lemma~\ref{lem:derivatives} are recovered when the corresponding eigenspaces are one-dimensional.

\begin{lemma}[Bounds on the second-order statistic $\xi$]\label{lem:ingredients}
The quantities in~\Thm{thm:RAP} satisfy $\mu^2\leq\xi\leq\mu$.
\end{lemma}
\begin{proof}
By Lemma~\ref{lem:loewner}, $\C\preccurlyeq\B/2$, and therefore
\begin{equation*}
    \xi=\langle\mU,\C(\mU)\rangle\leq\frac12\langle\mU,\B(\mU)\rangle=\frac12\lambda_1(\B)=\mu.
\end{equation*}
On the other hand, since $\xi=\expect{\left(\vect{u}^T\proj\vect{u}\right)^2}$ and $\mu=\expect{\vect{u}^T\proj\vect{u}}$, Jensen's inequality yields $\xi\geq\mu^2$.
\end{proof}

\subsection{The A-bound}\label{sec:a_bound}

The Loewner inequality $\C\preccurlyeq\B/2$ yields
the classical B-bound, but obtaining a sharper bound
requires exploiting the additional second-order
information contained in $\xi$.
Rather than seeking a surrogate that dominates $\C$
in the Loewner order, we seek a rank-one surrogate
$\C^\star$ such that the minimum Rayleigh quotient
of $\B-\omega\C^\star$ provides a lower bound
for that of $\B-\omega\C$.

We recall the structural properties of $\C$
established in the preceding subsections:
$\C\succcurlyeq\mat{0}$,
$\C\preccurlyeq\B/2$, and
$(\B-2\C)(\I)=0$, so that the quadratic forms
of $\C$ and $\B/2$ touch along the identity
direction. Furthermore,
$\langle\mU,\C(\mU)\rangle=\xi$.
These properties form the basis for constructing
and comparing the two-dimensional surrogates.

\begin{figure}
    \centering
    \includegraphics[width=0.5\linewidth]
        {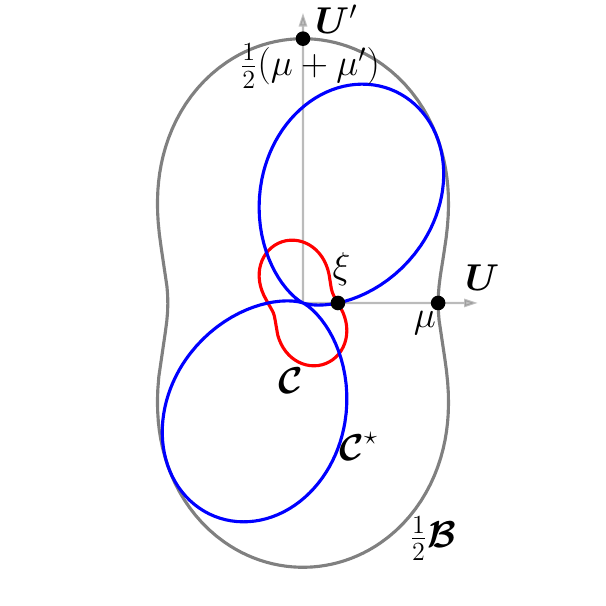}
    \caption{Rayleigh quotients of $\B/2$, $\C$,
    and $\C^\star$ restricted to the distinguished
    two-dimensional subspace
    $S^\star=\operatorname{span}\{\mU^\prime,\mU\}$.
    The surrogate $\C^\star$ need not dominate
    the compression of $\C$ in the Loewner order.}
    \label{fig:eclipse_c_star}
\end{figure}

We construct the surrogate $\C^\star$ using the
distinguished two-dimensional subspace
$S^\star=\operatorname{span}\{\mU^\prime,\mU\}$
identified in Fact~\ref{prop:BC}. In the ordered
orthonormal basis $(\mU^\prime,\mU)$, the compression
of $\B$ to $S^\star$ is represented by the
$2\times2$ matrix $\mat{B}$ in~\Thm{thm:RAP}.
The only additional information taken from $\C$
is the second-order statistic
$\xi=\langle\mU,\C(\mU)\rangle
=\expect{(\vect{u}^T\proj\vect{u})^2}$.
Using these quantities, we define the rank-one
surrogate
\begin{equation}\label{eq:c_star}
    \C^\star :=
    \left(\sqrt{\gamma}\mU^\prime+\sqrt{\xi}\mU\right)
    \otimes
    \left(\sqrt{\gamma}\mU^\prime+\sqrt{\xi}\mU\right),
    \qquad
    \gamma=\frac{\mu+\mu^\prime}{2}
    \left(1-\frac{\xi}{\mu}\right).
\end{equation}
The surrogate $\C^\star$ is supported on $S^\star$,
and its compression to this subspace is represented
by the $2\times2$ matrix $\mat{C}$ in~\Thm{thm:RAP}.
Figure~\ref{fig:eclipse_c_star} illustrates the
corresponding quadratic forms on this subspace.
Unlike a comparison in the Loewner order, our
approach does not require $\C^\star$ to dominate
$\C$, even when compressed to $S^\star$.

The difficulty in establishing the bound is that the
minimizing eigenvector $\mat{V}(\omega)$ of
$\B-\omega\C$ generally lies outside the distinguished
subspace $S^\star$. In the irreducible setting,
Perron--Frobenius theory implies that this eigenvector
is positive definite, and hence full rank, for
$\omega\in(0,2)$. We therefore need to extend the
comparison from $S^\star$ to the full space.

To this end, we consider arbitrary two-dimensional
subspaces containing $\mU$. For any unit matrix
$\mat{V}$ orthogonal to $\mU$, let
$S_{\mat{V}}:=\operatorname{span}\{\mat{V},\mU\}$,
and denote the compression of $\B$ to this subspace
by $\B_{\mat{V}}:=\B_{S_{\mat{V}}}$.
Since $\mU$ is an eigenvector of $\B$ corresponding
to $2\mu$, this compression is diagonal in the
ordered orthonormal basis $(\mat{V},\mU)$.

Writing $\pi_{\mat{V}}:=
\frac12\langle\mat{V},\B(\mat{V})\rangle$,
we construct the corresponding rank-one surrogate
\begin{equation}\label{eq:c_star_v}
    \C^\star_{\mat{V}}
    :=\left(\sqrt{\gamma_{\mat{V}}}\mat{V}
    +\sqrt{\xi}\mU\right)
    \otimes
    \left(\sqrt{\gamma_{\mat{V}}}\mat{V}
    +\sqrt{\xi}\mU\right),
    \qquad
    \gamma_{\mat{V}}
    :=\pi_{\mat{V}}\left(1-\frac{\xi}{\mu}\right).
\end{equation}
The parameter $\gamma_{\mat{V}}$ is nonnegative
by Lemma~\ref{lem:ingredients}. In the ordered
basis $(\mat{V},\mU)$, the compressed operators
are represented by
\[
    \mat{B}_{\mat{V}}=
    \begin{bmatrix}
        2\pi_{\mat{V}} & 0\\
        0 & 2\mu
    \end{bmatrix},
    \qquad
    \mat{C}^\star_{\mat{V}}=
    \begin{bmatrix}
        \gamma_{\mat{V}} &
        \sqrt{\gamma_{\mat{V}}\xi}\\
        \sqrt{\gamma_{\mat{V}}\xi} & \xi
    \end{bmatrix}.
\]
The surrogate $\C^\star_{\mat{V}}$ is supported
on $S_{\mat{V}}$, and choosing
$\mat{V}=\mU^\prime$ recovers the distinguished
surrogate $\C^\star$ defined in~\eqref{eq:c_star}.

\begin{proposition}[Properties of the rank-one surrogate]
\label{pro:c_star_v_kissing}
For any unit matrix $\mat{V}\perp\mU$, the surrogate
$\C^\star_{\mat{V}}$ is positive semidefinite and satisfies
$\langle\mU,\C^\star_{\mat{V}}(\mU)\rangle=\xi$.
Moreover, on $S_{\mat{V}}$, we have
$\C^\star_{\mat{V}}\preccurlyeq\frac12\B_{\mat{V}}$,
with $\frac12\B_{\mat{V}}-\C^\star_{\mat{V}}$ singular.
\end{proposition}

\begin{proof}
Positive semidefiniteness of $\C^\star_{\mat{V}}$
follows from its rank-one construction.
Since $\mat{V}\perp\mU$, we also have
$\langle\mU,\C^\star_{\mat{V}}(\mU)\rangle=\xi$.

In the ordered basis $(\mat{V},\mU)$, we have
\[
\frac12\mat{B}_{\mat{V}}-\mat{C}^\star_{\mat{V}}
=
\begin{bmatrix}
\pi_{\mat{V}}-\gamma_{\mat{V}} &
-\sqrt{\gamma_{\mat{V}}\xi}\\
-\sqrt{\gamma_{\mat{V}}\xi} & \mu-\xi
\end{bmatrix}.
\]
Both diagonal entries are nonnegative, since
$\xi\leq\mu$ (due to Lemma~\ref{lem:ingredients}) and
$\gamma_{\mat{V}}=\pi_{\mat{V}}(1-\xi/\mu)$.
Moreover, its determinant is
\[
(\pi_{\mat{V}}-\gamma_{\mat{V}})(\mu-\xi)
-\gamma_{\mat{V}}\xi
=\pi_{\mat{V}}(\mu-\xi)-\gamma_{\mat{V}}\mu=0.
\]
Thus, the matrix is positive semidefinite and singular,
establishing
$\C^\star_{\mat{V}}\preccurlyeq\frac12\B_{\mat{V}}$.
\end{proof}

The compression of $\C$ to $S_{\mat{V}}$ inherits
the properties $\C\succcurlyeq\mat{0}$,
$\C\preccurlyeq\B/2$, and
$\langle\mU,\C(\mU)\rangle=\xi$.
The preceding proposition shows that the rank-one
surrogate $\C^\star_{\mat{V}}$ satisfies these same
constraints on $S_{\mat{V}}$.
We now show that this surrogate yields a minimum
Rayleigh quotient no greater than that of the
compression of $\B-\omega\C$, uniformly for
$\omega\in[0,2]$.

\begin{proposition}[Comparison within a two-dimensional subspace]
\label{pro:c_star_v}
For any unit matrix $\mat{V}\perp\mU$, the rank-one
surrogate $\C^\star_{\mat{V}}$ satisfies
\begin{equation*}
\lambda_{\min}(\B_{\mat{V}}-\omega\C^\star_{\mat{V}})
\leq
\lambda_{\min}(\B_{\mat{V}}-\omega\C_{S_{\mat{V}}})
\end{equation*}
for every $\omega\in[0,2]$.
\end{proposition}

\begin{proof}
Since $\mat{0}\preccurlyeq\C\preccurlyeq\frac12\B$, its compression satisfies
$0\preccurlyeq\C_{S_{\mat{V}}}
\preccurlyeq\frac12\B_{\mat{V}}$.
However, unlike the rank-one surrogate in
Proposition~\ref{pro:c_star_v_kissing},
this compression need not touch $\frac12\B_{\mat{V}}$.
To obtain the same structural condition, we enlarge
the compression along $\mat{V}$, defining
$\widehat{\C}_{\mat{V}}
:=\C_{S_{\mat{V}}}+\delta\mat{V}\otimes\mat{V}$,
where $\delta\geq0$ is chosen as large as possible
while preserving
$\widehat{\C}_{\mat{V}}
\preccurlyeq\frac12\B_{\mat{V}}$.
Consequently, $\frac12\B_{\mat{V}}-\widehat{\C}_{\mat{V}}$
is positive semidefinite and singular. Moreover,
the enlargement preserves
$\langle\mU,\widehat{\C}_{\mat{V}}(\mU)\rangle=\xi$.

In the ordered basis $(\mat{V},\mU)$, the matrix
representation of $\widehat{\C}_{\mat{V}}$ can
therefore be written, up to the sign of its
off-diagonal entries, as
\[
\widehat{\mat{C}}_{\mat{V}}
=
\begin{bmatrix}
\alpha+\sigma & \sqrt{\alpha\xi}\\
\sqrt{\alpha\xi} & \xi
\end{bmatrix},
\qquad
\alpha=(\pi_{\mat{V}}-\sigma)
\left(1-\frac{\xi}{\mu}\right),
\quad 0\leq\sigma\leq\pi_{\mat{V}}.
\]
The choice $\sigma=0$ recovers the rank-one
surrogate $\C^\star_{\mat{V}}$.
As established by the monotonicity calculation
in Appendix~\ref{app:eclipse-proofs}, the smallest
eigenvalue of
$\B_{\mat{V}}-\omega\widehat{\C}_{\mat{V}}$
is nondecreasing in $\sigma$ for every
$\omega\in[0,2]$. Hence
\[
\lambda_{\min}(\B_{\mat{V}}-\omega\C^\star_{\mat{V}})
\leq
\lambda_{\min}(\B_{\mat{V}}-\omega\widehat{\C}_{\mat{V}})
\leq
\lambda_{\min}(\B_{\mat{V}}-\omega\C_{S_{\mat{V}}}),
\]
where the second inequality follows from
$\widehat{\C}_{\mat{V}}\succcurlyeq
\C_{S_{\mat{V}}}$.
\end{proof}

Having identified the optimal rank-one surrogate
within each two-dimensional subspace, we now compare
the choices of subspace. These surrogates depend on
$\mat{V}$ through
$\pi_{\mat{V}}:=\frac12\langle\mat{V},\B(\mat{V})\rangle$.
By the min--max theorem,
\[
\pi_{\mat{V}}\geq\frac12\lambda_2(\B)
=\frac{\mu+\mu^\prime}{2}.
\]
The minimum is attained when $\mat{V}=\mU^\prime$,
corresponding to the distinguished subspace $S^\star$.

\begin{proposition}[Comparison across two-dimensional subspaces]
\label{pro:c_star}
For any unit matrix $\mat{V}\perp\mU$ and every
$\omega\in[0,2]$, we have
\begin{equation*}
\lambda_{\min}(\B_{S^\star}-\omega\C^\star_{S^\star})
\leq
\lambda_{\min}(\B_{\mat{V}}-\omega\C^\star_{\mat{V}}).
\end{equation*}
\end{proposition}

\begin{proof}
As established in Appendix~\ref{app:eclipse-proofs},
the smallest eigenvalue of
$\B_{\mat{V}}-\omega\C^\star_{\mat{V}}$
is nondecreasing in $\pi_{\mat{V}}$ for every
$\omega\in[0,2]$. Since $\pi_{\mat{V}}$ is minimized
at $\mat{V}=\mU^\prime$, for which
$\C^\star_{\mat{V}}=\C^\star$, the result follows.
\end{proof}

The preceding propositions establish that, for every
two-dimensional subspace containing $\mU$, the corresponding
rank-one surrogate provides a lower bound on the minimum
Rayleigh quotient of the compressed operator, and that
the distinguished subspace $S^\star$ yields the strongest
such bound. Both comparisons hold for all $\omega\in[0,2]$.
We now extend these comparisons to the full space by
selecting, for each fixed $\omega$, a two-dimensional
subspace containing a minimizing eigenvector of
$\B-\omega\C$.

\begin{theorem}\label{thm:eclipse}
The rank-one surrogate $\C^\star$ defined
in~\eqref{eq:c_star} satisfies
\begin{equation}
\lambda_{\min}(\B-\omega\C^\star)
\leq
\lambda_{\min}(\B-\omega\C)
\end{equation}
for every $\omega\in[0,2]$.
\end{theorem}

\begin{proof}
Fix $\omega\in[0,2]$ and let $\mat{X}$ be a
normalized eigenvector corresponding to
$\lambda_{\min}(\B-\omega\C)$.

If $\mat{X}$ is proportional to $\mU$, then
$\lambda_{\min}(\B-\omega\C)=2\mu-\omega\xi$.
Since $\langle\mU,\C^\star(\mU)\rangle=\xi$,
the variational characterization gives
$\lambda_{\min}(\B-\omega\C^\star)
\leq 2\mu-\omega\xi$, establishing the result
in this case.

Otherwise, let $\mat{V}$ be the normalized
component of $\mat{X}$ orthogonal to $\mU$,
and consider the two-dimensional subspace
$S_{\mat{V}}=\operatorname{span}\{\mat{V},\mU\}$.
Since $\mat{X}\in S_{\mat{V}}$, the variational
characterization in~\eqref{eq:variational_slice} gives
\[
\lambda_{\min}(\B_{\mat{V}}
-\omega\C_{S_{\mat{V}}})
=
\lambda_{\min}(\B-\omega\C).
\]

Propositions~\ref{pro:c_star} and~\ref{pro:c_star_v}
establish
\begin{align*}
\lambda_{\min}(\B_{S^\star}
-\omega\C^\star_{S^\star})
&\leq
\lambda_{\min}(\B_{\mat{V}}
-\omega\C^\star_{\mat{V}})\\
&\leq
\lambda_{\min}(\B_{\mat{V}}
-\omega\C_{S_{\mat{V}}})\\
&=
\lambda_{\min}(\B-\omega\C).
\end{align*}

Since $S^\star$ is spanned by eigenvectors of $\B$,
it is an invariant subspace of $\B$, as is its
orthogonal complement. Moreover, $\C^\star$ is
supported on $S^\star$, so $\B-\omega\C^\star$
agrees with $\B$ on $(S^\star)^\perp$.
The smallest eigenvalue on this complement is
at least $2\mu$, whereas the minimum Rayleigh
quotient on $S^\star$ is at most
$2\mu-\omega\xi\leq 2\mu$.
Consequently, the smallest eigenvalue of
$\B-\omega\C^\star$ is attained on $S^\star$, and
\[
\lambda_{\min}(\B-\omega\C^\star)
=
\lambda_{\min}(\B_{S^\star}
-\omega\C^\star_{S^\star}).
\]
Combining these relations completes the proof.
\end{proof}

\begin{proof}[Proof of Theorem~\ref{thm:RAP}]
Since $\C^\star$ is supported on the subspace
$S^\star=\operatorname{span}\{\mU^\prime,\mU\}$,
the smallest eigenvalue of $\B-\omega\C^\star$
is given by the $2\times2$ matrices $\mat{B}$ and
$\mat{C}$ in~\Thm{thm:RAP}. Therefore,
Theorem~\ref{thm:eclipse} establishes
\[
\lambda_{\min}(\mat{B}-\omega\mat{C})
=\lambda_{\min}(\B-\omega\C^\star)
\leq\lambda_{\min}(\B-\omega\C).
\]
Combining this inequality with~\eqref{eq:AtoBC}
and Proposition~\ref{prop:lyapunov_jensen} yields
\begin{equation}
\varphi(\omega)\leq\overline{\varphi}(\omega)
\leq\AB(\omega).
\end{equation}
\end{proof}

Furthermore, since $\C^\star\preccurlyeq\frac12\B$,
we have $\AB(\omega)\leq\BB(\omega)$.
The strength of our approach, compared with bounds
based on the Loewner order, is that it does not
require $\C^\star$ to dominate $\C$.
Indeed, as illustrated in~\FG{fig:eclipse_c_star},
such domination need not hold even when both
superoperators are compressed to the distinguished
subspace $S^\star$.

\section{Applications in Iterative Methods}
\label{sec:app}

Classical stationary iterative methods for solving linear systems $\mA\bx=\bb$ have long been studied through their asymptotic convergence factors \cite{saad2003iterative,hackbusch2016iterative,varga1962iterative}. For a stationary method with iteration matrix $\mT$, the asymptotic convergence factor is given by its spectral radius, $\rho(\mT)=\lim_{k\to\infty}\|\mT^k\|^{1/k}$. Relating this spectral radius to properties of $\mA$ is generally difficult because the iteration matrix depends on the algorithm, including the splitting and ordering of variables. Even for classical methods such as Gauss--Seidel and successive over-relaxation (SOR), explicit relationships are available only under special structural assumptions~\cite{berman1994nonnegative}. Similarly, for deterministic Kaczmarz and other alternating-projection methods, convergence rates depend on the angles between subspaces and the order in which projections are applied. Bounds based on these angles are often difficult to compute and relate to spectral properties of $\mA$, such as its condition number~\cite{galantai2005rate,deutsch1997rate}.

Randomization does not, by itself, resolve the difficulty of characterizing asymptotic convergence factors. Replacing a fixed iteration matrix by a random product leads to the largest Lyapunov exponent, which is notoriously difficult to compute or approximate in general~\cite{tsitsiklis1997lyapunov,blondel2005computationally,altschuler2020lyapunov}. Nevertheless, randomization enables a tractable per-iteration analysis, yielding computable convergence guarantees for methods such as randomized Kaczmarz and Gauss--Seidel~\cite{strohmer2009randomized,leventhal2010randomized}. As discussed in~\SC{sec:intro}, these guarantees can be considerably more conservative than the observed asymptotic factors and fail to predict improvements through overrelaxation. We now apply~\Thm{thm:RAP} to randomized Kaczmarz, Gauss--Seidel, gossip algorithms, and block subspace corrections, quantifying improvements in convergence guarantees both without relaxation and through overrelaxation.

\subsection*{Iteration Complexity}

The practical significance of a sharper convergence factor lies in the number of iterations sufficient to guarantee a prescribed accuracy~\cite{nesterov2012efficiency}. The spectral radius of the second-moment superoperator provides the following finite-iteration guarantee.

\begin{proposition}[Finite-iteration convergence guarantee]
\label{prop:finite-iteration}
Under the assumptions of Theorem~\ref{thm:RAP}, for any initial error $\mathbf{z}_0\in\mathbb{R}^n$, any $\omega\in(0,2)$, and every $k\geq0$, the expected squared error satisfies
\begin{equation}
\mathbb{E}\|\mathbf{z}_k\|^2
\leq \sqrt{n}\,\overline{\varphi}(\omega)^k
\|\mathbf{z}_0\|^2
\leq \sqrt{n}\,\AB(\omega)^k
\|\mathbf{z}_0\|^2.
\end{equation}
Consequently, for a relative expected squared-error tolerance $\varepsilon\in(0,1)$, the A-bound guarantees the prescribed accuracy after
\begin{equation}
k_{\mathrm A}(\varepsilon,\omega)
:=
\left\lceil
\frac{\log(\sqrt{n}/\varepsilon)}
{-\log\AB(\omega)}
\right\rceil
\end{equation}
iterations, provided $0<\AB(\omega)<1$.
\end{proposition}

\begin{proof}
Recall that $\mat{\Sigma}_k=\A^k(\mat{\Sigma}_0)$, where $\mat{\Sigma}_0=\mathbf{z}_0\mathbf{z}_0^T$. Since $\mat{\Sigma}_k$ is positive semidefinite,
\begin{align*}
\mathbb{E}\|\mathbf{z}_k\|^2
&=\operatorname{tr}(\mat{\Sigma}_k)
\leq \sqrt{n}\|\mat{\Sigma}_k\|_F
\leq \sqrt{n}\|\A^k\|\|\mat{\Sigma}_0\|_F
=\sqrt{n}\rho(\A)^k\|\mathbf{z}_0\|^2\\
&\leq \sqrt{n}\AB(\omega)^k\|\mathbf{z}_0\|^2.
\end{align*}
Here, the operator norm of $\A$ is induced by the Frobenius norm and coincides with its spectral radius due to self-adjointness. The final inequality follows from Theorem~\ref{thm:RAP}. Solving for $k$ establishes the iteration-complexity guarantee.
\end{proof}

For comparison, the classical B-bound at $\omega=1$ guarantees the same tolerance after
\[
k_{\mathrm B}(\varepsilon)
:=
\left\lceil
\frac{\log(1/\varepsilon)}
{-\log\BB(1)}
\right\rceil
\]
iterations. Let $\omega_\star$ denote the relaxation parameter minimizing $\AB(\omega)$ over $\omega\in(0,2)$. We distinguish the unrelaxed guarantee $k_{\mathrm A}(\varepsilon,1)$ from the optimally relaxed guarantee $k_{\mathrm A}(\varepsilon,\omega_\star)$.

In the examples below, we quantify the improvements in convergence guarantees in three ways: the improvement in guaranteed contraction provided by the A-bound over the classical B-bound at $\omega=1$, the additional improvement achieved through A-bound-optimal overrelaxation, and the resulting reduction in the number of iterations sufficient to guarantee a prescribed accuracy. Throughout the numerical examples in this section, we fix the relative expected squared-error tolerance at $\varepsilon=10^{-8}$ when reporting iteration-complexity guarantees.

\subsection{Randomized Kaczmarz}\label{sec:kacz}

For a consistent system of $m$ linear equations $\mA\bx=\bb$, randomized Kaczmarz updates the iterate by selecting an equation and applying a relaxed projection onto its associated affine hyperplane. Let $\vect{a}_i:=\mA^T\be_i$, where $\be_i$ denotes the $i$th canonical basis vector. At iteration $k$, an index $i_k$ is selected randomly, and the iterate is updated according to
\begin{equation}
\bx_k=\bx_{k-1}+\omega\frac{\be_{i_k}^T(\bb-\mA\bx_{k-1})}{\vect{a}_{i_k}^T\vect{a}_{i_k}}\vect{a}_{i_k}.
\end{equation}
The residual $\be_{i_k}^T(\bb-\mA\bx_{k-1})$ can be computed in $O(n)$ floating-point operations, and $\omega\in(0,2)$ is the relaxation parameter. Defining the orthogonal projectors $\proj_i:=\vect{a}_i\vect{a}_i^T/(\vect{a}_i^T\vect{a}_i)$ and the error vector $\vect{z}_k:=\bx_k-\bx^\star$, we recover the random product of projections formulation~(\ref{eq:factor}): $\vect{z}_k=(I-\omega\proj_{i_k})\vect{z}_{k-1}.$

The expected projector can be directly related to the underlying matrix $\mA$. Selecting row $i$ with probability $p_i:=\|\vect{a}_i\|^2/\|\mA\|_F^2$ yields $\mathbb{E}[\proj]=\mA^T\mA/\|\mA\|_F^2$. Thus, the classical bound $\BB(1)=1-\mu$ is governed by the scaled condition number $\kappa:=1/\mu=\|\mA\|_F^2/\sigma_{\min}^2(\mA)$. When the rows of $\mA$ have unit norm, this distribution reduces to uniform sampling. In this case, the three ingredients of the A-bound in~\Thm{thm:RAP} admit simple expressions in terms of the singular values and singular vectors of $\mA$: $\mu$ and $\mu'$ are the two smallest squared singular values of $\mA/\|\mA\|_F$, and $\xi=m\mu^2\|\vect{v}\|_4^4$, where $\vect{v}$ is the unit left singular vector corresponding to the smallest singular value of $\mA$.

\begin{figure}[ht!]
\centering
\includegraphics[width=0.96\textwidth]{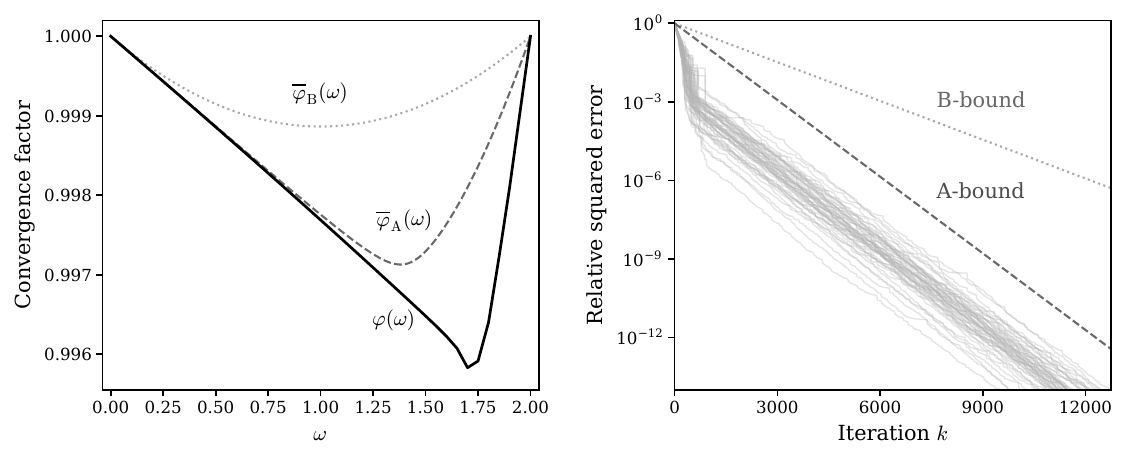}
\caption{Convergence of randomized Kaczmarz with uniform sampling on the nonsymmetric Parter matrix for $n=100$ $(\kappa=880.3)$. The left panel compares the asymptotic convergence factor $\varphi(\omega)$ with the A-bound $\AB(\omega)$ and the classical B-bound $\BB(\omega)$. The right panel shows the error decay from 50 independent runs with $\omega=1$ on a logarithmic scale, together with reference lines corresponding to $\AB(1)$ and $\BB(1)=1-\mu$. The slope of the A-bound substantially narrows the gap between the classical guarantee and the observed asymptotic convergence factor.}
\label{fig:Parter}
\end{figure}

\Cref{fig:Parter} shows the performance of randomized Kaczmarz with uniform sampling on the nonsymmetric Parter matrix~\cite{parter1986distribution,moler} for $n=100$. The left panel compares the convergence bounds with the observed asymptotic factor, while the right panel shows the error decay on a logarithmic scale over an ensemble of 50 independent trials. After an initial rapid decrease, the error trajectories settle into nearly parallel linear decay, consistent with a common asymptotic convergence factor $\varphi(1)$. At $\omega=1$, the A-bound substantially narrows the gap between the classical B-bound and the observed asymptotic convergence factor, improving the guaranteed contraction by a factor of approximately $1.98$. The predicted optimal relaxation parameter $\omega_\star\approx1.38$ yields a further $28\%$ improvement in guaranteed contraction, or an overall improvement by a factor of approximately $2.53$ relative to the classical B-bound. The corresponding iteration-complexity guarantees are $k_{\mathrm B}(\varepsilon)\approx16{,}207$, $k_{\mathrm A}(\varepsilon,1)\approx9{,}225$, and $k_{\mathrm A}(\varepsilon,\omega_\star)\approx7{,}206$, representing reductions of $43.1\%$ and $55.5\%$, respectively, relative to the classical guarantee.

We next consider Gaussian matrices, a common benchmark in randomized Kaczmarz, for which the asymptotic bounds admit further analytical insights. For a matrix $\mA\in\mathbb{R}^{m\times n}$ with independent standard Gaussian entries and row-norm-proportional sampling, $\EP=\mA^T\mA/\|\mA\|_F^2$ is a normalized Wishart matrix. By rotational invariance, the unit left singular vector $\vect{v}$ associated with the smallest singular value satisfies $\mathbb{E}\|\vect{v}\|_4^4=3/(m+2)$, and concentration of Gaussian row norms yields $\xi/\mu^2\to3$ in probability as $m,n\to\infty$ with a fixed aspect ratio. The Marchenko--Pastur spectral behavior reveals an interesting distinction: for square Gaussian matrices, the hard-edge fluctuations of the smallest eigenvalues allow a significant separation between $\mu$ and $\mu'$, potentially leading to substantial benefits from overrelaxation. For a $500\times500$ realization, the effect of overrelaxation is particularly pronounced. At $\omega=1$, the A-bound improves the guaranteed contraction by a factor of approximately $2$ over the classical B-bound. The predicted optimal relaxation parameter $\omega_\star\approx1.77$ yields a further $76.5\%$ improvement, for an overall improvement by a factor of approximately $3.53$. The corresponding iteration-complexity guarantees are $k_{\mathrm B}(\varepsilon)\approx5.55\times10^{10}$, $k_{\mathrm A}(\varepsilon,1)\approx3.24\times10^{10}$, and $k_{\mathrm A}(\varepsilon,\omega_\star)\approx1.84\times10^{10}$, representing reductions of $41.6\%$ and $66.9\%$, respectively, relative to the classical guarantee. In contrast, for tall Gaussian matrices with a fixed aspect ratio $m/n>1$, the ratio $\mu'/\mu$ approaches one, and the additional improvement predicted by the A-bound vanishes asymptotically. This contrast is illustrated in \FG{fig:gaussian-gap}, which compares the asymptotic convergence factor and its bounds for $50\times50$ and $75\times50$ Gaussian matrices. Capturing potential improvements in this regime requires extending the analysis to a cluster of small eigenvalues, as discussed in~\SC{sec:discussion}.

\begin{figure}[t]
    \centering
    \includegraphics[width=\linewidth]{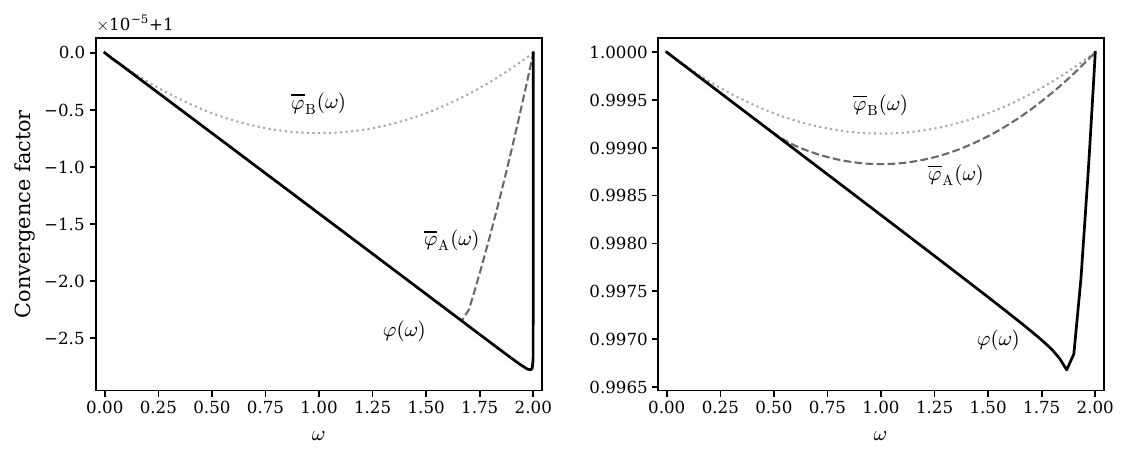}
    \caption{Randomized Kaczmarz examples illustrating how the spectral gap between $\mu$ and $\mu'$ influences the improvement of the A-bound over the B-bound. The left and right panels correspond to $50\times50$ and $75\times50$ Gaussian matrices, respectively, uniform sampling of normalized rows. The solid curves show the estimated asymptotic convergence factor $\varphi(\omega)$, while the dashed and dotted curves show the A- and B-bounds, respectively. The corresponding ratios are $\mu'/\mu\approx11.52$ (left) and $1.76$ (right).}
    \label{fig:gaussian-gap}
\end{figure}

As a third example, we examine randomized Kaczmarz in pairwise gossip algorithms for distributed averaging over networks~\cite{boyd2006randomized,loizou2016new}. For a connected, undirected graph with $n$ vertices and $m$ edges, let $\mA$ be its oriented edge--vertex incidence matrix and $\mL=\mA^T\mA$ the graph Laplacian. Uniform edge sampling yields $\EP=\mL/(2m)$, and convergence is analyzed on the $(n-1)$-dimensional orthogonal complement of the constant vector, where $\mu$ and $\mu'$ correspond to the Fiedler eigenvalue and the next eigenvalue of $\mL$, respectively, scaled by $1/(2m)$. For two five-vertex cliques connected by a single edge, the A-bound at $\omega=1$ improves the guaranteed contraction by a factor of approximately $1.83$, while the predicted optimal relaxation parameter $\omega_\star\approx1.46$ yields a further $23.5\%$ improvement. The corresponding iteration-complexity guarantees are $k_{\mathrm B}(\varepsilon)=2584$, $k_{\mathrm A}(\varepsilon,1)=1490$, and $k_{\mathrm A}(\varepsilon,\omega_\star)=1204$, representing reductions of $42.3\%$ and $53.4\%$, respectively, relative to the classical guarantee. For the dolphin social network~\cite{lusseau2003emergent}, the A-bound improves the guaranteed contraction at $\omega=1$ by a factor of $1.96$, but optimal overrelaxation at $\omega_\star\approx1.12$ provides a more modest additional improvement of $5.3\%$. The corresponding iteration counts are $33{,}856$, $19{,}205$, and $18{,}243$, respectively, yielding reductions of $43.3\%$ and $46.1\%$ relative to the classical guarantee.

\subsection{Randomized Gauss-Seidel}\label{sec:gs}
For a symmetric positive definite linear system $\mA\bx=\bb$, randomized Gauss--Seidel updates one coordinate of the iterate at each step. At iteration $k$, a coordinate $i_k$ is selected randomly, and the iterate is updated according to
\begin{equation}
\bx_k=\bx_{k-1}+\omega\frac{\be_{i_k}^T(\bb-\mA\bx_{k-1})}{\mA_{i_k i_k}}\be_{i_k},
\label{eq:gs_iter}
\end{equation}
where $\omega\in(0,2)$ is the relaxation parameter. The residual $\be_{i_k}^T(\bb-\mA\bx_{k-1})$ can be computed in $O(n)$ floating-point operations. 
The Gauss--Seidel iteration can be expressed in terms of the projectors $\widetilde{\proj}_i:=\be_i\be_i^T\mA/\mA_{ii}$, which are generally not orthogonal in the Euclidean inner product. When $\mA$ is symmetric positive definite, the change of coordinates $\proj_i:=\mA^{1/2}\widetilde{\proj}_i\mA^{-1/2}=\mA^{1/2}\be_i\be_i^T\mA^{1/2}/\mA_{ii}$ yields orthogonal projectors. Defining the transformed error $\vect{z}_k:=\mA^{1/2}(\bx_k-\bx^\star)$, we recover the random projection formulation~(\ref{eq:factor}): $\vect{z}_k=(\I-\omega\proj_{i_k})\vect{z}_{k-1}.$ In particular, $\|\vect{z}_k\|=\|\bx_k-\bx^\star\|_{\mA}$, so the convergence guarantees apply to the error in the $\mA$-norm.

The expected projector can be directly related to the underlying matrix $\mA$. Selecting coordinate $i$ with probability $p_i:=\mA_{ii}/\tr{\mA}$ yields $\EP=\mA/\tr{\mA}$. Consequently, $\mu$ and $\mu'$ are the two smallest eigenvalues of $\mA/\tr{\mA}$. When $\mA$ has unit diagonal, this distribution reduces to uniform sampling, and $\xi=n\mu^2\|\vect{u}\|_4^4$, where $\vect{u}$ is a unit eigenvector associated with the smallest eigenvalue of $\mA$.

\begin{figure}[ht!]
\centering
\begin{tabular}{@{\hspace{0pt}}c@{\hspace{1pt}}c}
\myfig{.45}{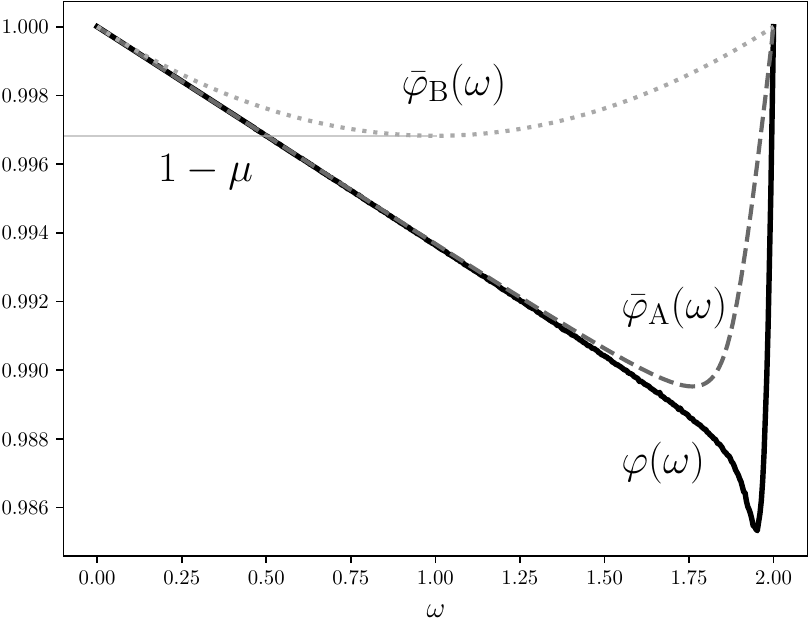}&
\myfig{.45}{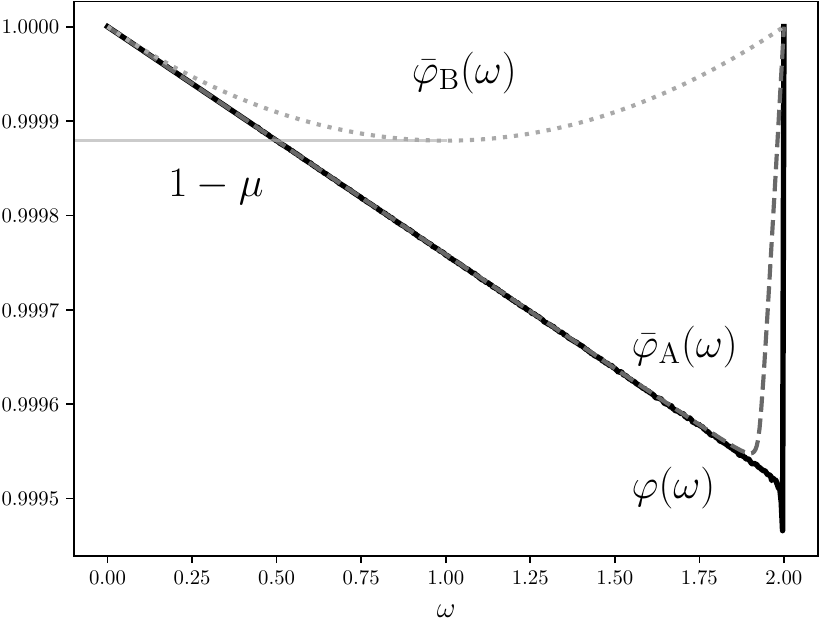}
\end{tabular}
\caption{Convergence of randomized Gauss--Seidel with uniform sampling on Hilbert matrices of order $n=3$ (left, $\kappa=314.7$) and $n=4$ (right, $\kappa=8280.5$). The asymptotic convergence factor $\varphi(\omega)$ is compared with the A-bound $\AB(\omega)$ and the classical B-bound $\BB(\omega)$. Additional experiments for $n>4$ (not shown) reveal a further narrowing of the gap between $\AB$ and $\varphi$, to the point where the two curves become nearly indistinguishable.}
\label{fig:Hilbert}
\end{figure}

\Cref{fig:Hilbert} illustrates how the A-bound becomes increasingly accurate as the conditioning of the Hilbert matrix worsens. The agreement extends to the predicted optimal relaxation parameter, with the minimizer of $\AB(\omega)$ approaching that of the true asymptotic convergence factor $\varphi(\omega)$. Additional experiments with $n>4$ (not shown) exhibit the same trend, with the A-bound and the true convergence factor becoming nearly indistinguishable.

As a larger-scale sparse SPD example, we consider the {\tt bcsstk09} matrix from the SuiteSparse Matrix Collection, a structural engineering stiffness matrix of order $n=1083$. We apply randomized Gauss--Seidel to the original matrix, selecting coordinate $i$ with probability $p_i=A_{ii}/\operatorname{tr}(A)$. At $\omega=1$, the A-bound improves the guaranteed contraction by a factor of approximately $2$ over the classical B-bound. The predicted optimal relaxation parameter $\omega_\star\approx1.18$ yields a further $17.4\%$ improvement, for an overall improvement by a factor of approximately $2.35$. The corresponding iteration-complexity guarantees are $k_{\mathrm B}(\varepsilon)\approx4.75\times10^7$, $k_{\mathrm A}(\varepsilon,1)\approx2.83\times10^7$, and $k_{\mathrm A}(\varepsilon,\omega_\star)\approx2.41\times10^7$, representing reductions of $40.5\%$ and $49.3\%$, respectively, relative to the classical guarantee.

\subsection{Randomized Block Subspace Iterations}\label{sec:block}

The applicability of~\Thm{thm:RAP} is not restricted to rank-one projectors. We now consider randomized block subspace corrections, where each iteration applies a correction over a potentially higher-dimensional subspace. This framework includes randomized block Gauss--Seidel and Schwarz-type methods arising in domain decomposition.

For symmetric positive definite systems, exact block Gauss--Seidel corrections can be represented as orthogonal projections in the $\mA$-geometry, with $$\proj_J=\mA^{1/2}\mE_J(\mE_J^T\mA \mE_J)^{-1}\mE_J^T\mA^{1/2},$$ where $\mE_J$ selects the coordinates in block $J$. This extends the rank-one formulation of the preceding subsection. Randomized block methods and their sampling strategies have been studied extensively~\cite{nesterov2012efficiency, tu2017breaking, gower2015randomized}. Since the resulting projectors may have rank greater than one, this setting provides a natural application of~\Thm{thm:RAP} beyond rank-one projections.

As a concrete example, consider the Poisson equation $-\Delta u=f$ on $(0,1)^2$, with Dirichlet boundary conditions on the left boundary and homogeneous Neumann conditions on the remaining boundaries. A standard five-point finite-difference discretization on a $12\times12$ grid yields an SPD system with $144$ unknowns. We partition the unknowns into nonoverlapping rectangular blocks and apply exact block corrections, selecting blocks uniformly at random.

With four blocks of $36$ unknowns each, the A-bound at $\omega=1$ improves the guaranteed contraction by a factor of approximately $1.96$ over the classical B-bound. The predicted optimal relaxation parameter $\omega_\star\approx1.36$ yields a further $24.1\%$ improvement, for an overall improvement by a factor of approximately $2.44$. The corresponding iteration-complexity guarantees are $k_{\mathrm B}(\varepsilon)=1755$, $k_{\mathrm A}(\varepsilon,1)=1011$, and $k_{\mathrm A}(\varepsilon,\omega_\star)=812$, representing reductions of $42.4\%$ and $53.7\%$, respectively.

With six blocks of $24$ unknowns each, the A-bound improves the guaranteed contraction by a factor of approximately $1.98$, with a further $29.9\%$ improvement through overrelaxation at $\omega_\star\approx1.39$. The corresponding iteration counts are $3752$, $2145$, and $1649$, respectively, yielding reductions of $42.8\%$ and $56.1\%$ relative to the classical guarantee.

We also revisit the {\tt bcsstk09} matrix from~\SC{sec:gs}, partitioning its $1083$ variables into six nonoverlapping blocks of approximately equal size, consisting of consecutive variables in the original ordering. With uniform block sampling, the A-bound at $\omega=1$ improves the guaranteed contraction by a factor of approximately $2$. The predicted optimal relaxation parameter $\omega_\star\approx1.20$ yields a further $17.5\%$ improvement, for an overall improvement by a factor of approximately $2.35$ relative to the classical B-bound. The corresponding iteration-complexity guarantees are $k_{\mathrm B}(\varepsilon)=108{,}625$, $k_{\mathrm A}(\varepsilon,1)=64{,}662$, and $k_{\mathrm A}(\varepsilon,\omega_\star)=55{,}028$, representing reductions of $40.5\%$ and $49.3\%$, respectively, relative to the classical guarantee.

These results demonstrate that the benefits of the asymptotic analysis extend beyond rank-one projections to higher-dimensional subspace corrections. These convergence guarantees are expressed per randomized block update; the computational cost of each update depends on the block size and the implementation of the local solves.

\section{Discussion and Related Work}\label{sec:discussion}
Stein and Rosenberg~\cite{stein1948solution} established the application of asymptotic analysis in the deterministic regime, which became foundational to the analysis of iterative methods~\cite{saad2003iterative}, supporting principled algorithm design and enabling comparison theorems~\cite{varga1962iterative} that allow one to formally determine when one method is faster than another (e.g., Gauss–Seidel vs. Jacobi). The A-bound introduced in~\Thm{thm:RAP} provides a corresponding asymptotic analysis in the randomized setting.

In the deterministic regime, the method of alternating projections is well studied~\cite{escalante2011alternating}. The fundamental difficulty in deriving the rate of convergence lies in the geometric nature of iterations. It is well known that the rate of convergence not only depends on angles between subspaces these projectors represent, but more importantly on the order these subspaces are visited~\cite{deutsch1995angle, galantai2005rate}. In this context a bound on the rate of convergence~\cite{halperin1962product} can be computed based on a sequence of angles between subspaces as they appear in the sweep~\cite{deutsch1997rate, ssw1977, kayalar1988error}. These bounds are difficult to calculate~\cite{galantai2005rate} and specific to the order in which subspaces are visited. 

Bauschke~\cite{bauschke1995norm} established norm convergence of products of relaxed projections in Hilbert spaces under suitable regularity conditions. However, convergence of alternating projections can be arbitrarily slow in infinite-dimensional settings~\cite{bauschke2009slow}.

The over-relaxation phenomenon, introduced by Young~\cite{young1954iterative} has played a major role in (deterministic) iterative methods. Interestingly Perron-Frobenius theory of (element-wise) positive matrices plays a foundational role in convergence analysis~\cite{varga1962iterative}. However, the question of optimal over-relaxation $\omega$ proves difficult and can only be answered for special matrices that satisfy certain structural conditions (e.g., consistently ordered cyclic matrices)~\cite{wild1987over, berman1994nonnegative, hackbusch2016iterative}. The randomization presents a compelling case as over-relaxations with improved and quantified factors, albeit not optimal, can be obtained in more general settings.

A common technique for bounding the spectral radius of Perron-Frobenius matrices is based on the sub-invariance theorem~\cite[Theorem~1.6]{seneta2006non}. For a positive operator (e.g., positive entries or positive linear maps) an approximate Perron eigenvector $\mat{V}\succ\mat{0}$ provides a bound on the spectral radius: $\A(\mat{V}) \preceq r \mat{V}$ implies $\rho(\A) \le r$ (see also~\cite{lessard2016analysis}). A better approximation to the Perron eigenvector produces a better bound on the radius. The A-bound developed in~\Thm{thm:RAP} is developed without guidance from the Perron eigenvector. It remains to be determined whether this bound corresponds to an approximate Perron eigenvector that can be characterized in closed-form.

The strength of the A-bound, compared to the B-bound,
depends on the gap between the eigenvalues $\mu$ and
$\mu'$. When the gap is small, the A-bound approaches
the B-bound. This motivates extending the analysis
to a cluster of small eigenvalues separated from the
remaining spectrum by a larger spectral gap.
Let $\mathcal{V}_{\mathrm{low}}$ denote the corresponding
spectral subspace of $\B$. A natural generalization of the second-order
statistic is
$\xi_{\mathrm{low}} :=
\max_{\mat{V}\in\mathcal{V}_{\mathrm{low}}\cap\mathbb{S}}
\langle\mat{V},\C(\mat{V})\rangle
=\lambda_{\max}(\C_{\mathcal{V}_{\mathrm{low}}})$.
This generalizes the statistic $\xi$ in
Theorem~\ref{thm:RAP}, recovered when
$\mathcal{V}_{\mathrm{low}}=\operatorname{span}\{\mU\}$. 
Repeated smallest eigenvalues constitute a special
case, in which $\mathcal{V}_{\mathrm{low}}$ is the
entire eigenspace of $\B$ corresponding to $2\mu$. For clusters of unequal eigenvalues, the largest
eigenvalue within the cluster determines its
spectral separation from the remaining spectrum,
while the variation of $\B$ within the cluster
must also be accounted for in constructing the bound.

In large-scale applications where the power of randomization becomes relevant\cite{strohmer2009randomized}, a realistic choice for $i_k$ is an ${\rm iid}$ process with uniform probability in $\{1, \dots, m\}$. Nevertheless, different sampling probabilities have been explored~\cite{nesterov2012efficiency, recht2012toward, agaskar2014randomized, gower2015randomized}. 

Asymptotic analysis of accelerated iterative methods that use an auxiliary variable $\vect{y}_k$ such as in accelerated Gauss-Seidel~\cite{nesterov2012efficiency, tu2017breaking} and accelerated Kaczmarz~\cite{lee2013efficient, liu2016accelerated} is of significant interest. Formulation of accelerated methods as nonsymmetric block matrices acting on $\bx_k$ and $\vect{y}_k$  provides a dynamical system view that can be analyzed in the asymptotic framework. The existing per-iteration analysis~\cite{nesterov2012efficiency, lee2013efficient, allen2016even, nesterov2017efficiency} provides a B-bound, $1 - \sqrt{{\mu}/{m}}$, that experimental evidence shows can be tightened in the asymptotic analysis similar to the unaccelerated iterations. Deriving asymptotic bounds, similar to the one provided in~\Thm{thm:RAP}, requires extending our min-max approach to the case of non-orthogonal eigenvectors that arise in the accelerated setting. Here the Birkhoff-Varga theorem provides a min-max formulation of the spectral radius for Perron-Frobenius operators in the nonsymmetric case. Nevertheless, a concrete A-bound in this case remains an open problem. An interesting question is whether the asymptotic analysis can lead to a better choice of parameters than  Nesterov's approach, since the benefits of overrelaxation manifest only in the asymptotic regime and remain invisible to the per-iteration analysis.

\appendix

\section{Perron-Frobenius Theory For Positive Linear Maps}\label{sec:PF}
The superoperator $\A$ defined in~\eqref{eq:covariance_dynamics} plays, for randomized iterations, the role of the iteration matrix in classical iterative methods, whose spectral radius governs asymptotic convergence~\cite{saad2003iterative}. In this section we discuss the theoretical foundations that provide necessary properties on the spectrum of $\A$ in the second-moment analysis we have seen in~\SC{sec:framework}.

Recall the superoperator $\A$, for a fixed $\omega$, denotes a linear map over the space of $n \times n$ matrices as:
\begin{equation*}
    \A(\mat{X}) = \expect{(\I - \omega \proj) \mat{X} (\I - \omega \proj)^T},
\end{equation*}
with the expectation taken over the set of $m$ projectors $\{\proj_i\}_{i=1}^m$ with probabilities $\{p_i\}_{i=1}^m$  for choosing indices $i_k$ in~\eqref{eq:factor}.
For any symmetric positive semi-definite matrix $\mat{X}$ the operation $(\I - \omega \proj_i) \mat{X} (\I - \omega \proj_i)^T$ preserves its positivity even when $\proj_i$ is not symmetric. Hence the superoperator $\A$ is a {\em positive linear map}~\cite{bhatia2009positive}, leaving the cone of symmetric positive semi-definite matrices invariant.

The spectra of positive linear maps on general (noncommutative) matrix algebras was studied in~\cite{evans1978spectral} that generalized the Perron-Frobenius theorem to this context. The spectral radius of a positive linear map is attained by an eigenvalue for which there exists an eigenvector that is positive semi-definite (see Theorem 6.5 in~\cite{wolf2012quantum}). The notion of {\em irreducibility} for positive linear maps guarantees that the eigenvalue is simple and the corresponding eigenvector is well-defined (up to a sign). What is more is that the eigenvector can be chosen to be a positive definite matrix. The corresponding positive-definite Perron eigenvector is associated with the spectral radius $\rho(\A)$, which provides the asymptotic bound on the second-moment dynamics $\mat{\Sigma}_k=\A^k(\mat{\Sigma}_0)$.

We examine the irreducibility of its corresponding superoperator $\A$ for any given relaxation value $\omega$. The criteria for irreducibility of positive linear maps was developed in~\cite{farenick1996irreducible} and involve invariant subspaces. A collection $S$ of (closed) subspaces of the vector space of $n \times n$ matrices is called trivial if it only contains $\{\mat{0}\}$ and the space itself. Given a bounded linear operator $\mat{M}$, let ${\rm Lat}(\mat{M})$ denote the invariant subspace lattice of $\mat{M}$. The following theorem is a specialization of a more general result in~\cite[Theorem 2]{farenick1996irreducible} to our superoperator.

\begin{theorem}[Irreducibility of the superoperator $\A$]
The positive linear map $\A$ is irreducible if and only if,  $\bigcap_{i=1}^m{{\rm Lat}(\I - \omega \proj_i})$ is trivial.
\end{theorem}
We note that the invariant subspaces of $(\I - \omega \proj_i)$ is the same as those of $\proj_i$.

Irreducibility is a convenient but not necessary assumption for the validity of the A-bound. Irreducibility guarantees that the spectral radius of $\A$ is unique for any relaxation value $\omega$. %
Irreducible problems  are dense in the sense that any such degenerate problem can be infinitesimally perturbed to a problem that the theorem applies to. Since the spectral radius of $\A$ is a continuous function of the underlying problem, the bounds remain valid regardless of irreducibility.

\section{Technical Proofs of the Two-Dimensional Comparisons}
\label{app:eclipse-proofs}

We provide the technical details for the two
comparisons established in~\SC{sec:a_bound}.
Both rely on monotonicity of the smallest eigenvalue
of a $2\times2$ symmetric matrix.
The first comparison varies $\sigma$ within
a fixed two-dimensional subspace, while the
second varies $\pi_{\mat{V}}$ across such subspaces.
Throughout, we use the two-dimensional matrix
representations associated with~\eqref{eq:c_star_v}.

\begin{proof}[Proof of Proposition~\ref{pro:c_star_v}]

Fix a unit matrix $\mat{V}$ orthogonal to $\mU$
and consider the subspace spanned by
$\mU$ and $\mat{V}$. Recall that $\pi_{\mat{V}}
:=\frac12\langle\mat{V},\B(\mat{V})\rangle$.

Let $\C_{\mat{V}}$ denote a positive semidefinite
operator on $S_{\mat{V}}$ satisfying
$\langle\mU,\C_{\mat{V}}(\mU)\rangle=\xi$,
such that $\frac12\B_{\mat{V}}-\C_{\mat{V}}$
is positive semidefinite and singular.
In particular, the saturated compression
$\widehat{\C}_{\mat{V}}$ constructed in the
main-text proof satisfies these properties.
Up to the sign of $\mat{V}$, such an operator
admits the following $2\times2$ representation
with $\alpha,\sigma\geq0$:
\begin{equation*}
    \mat{B}_{\mat{V}} = \begin{bmatrix}
        2\pi_{\mat{V}} & 0 \\ 0 & 2\mu
    \end{bmatrix} \quad {\rm and} \quad 
    \mat{C}_{\mat{V}} = \begin{bmatrix} \sqrt{\alpha} \\ \sqrt{\xi}\end{bmatrix} \begin{bmatrix} \sqrt{\alpha} & \sqrt{\xi}\end{bmatrix} + \begin{bmatrix} \sigma & 0 \\ 0 & 0 \end{bmatrix} = \begin{bmatrix} \alpha+\sigma & \sqrt{\alpha \xi} \\ \sqrt{\alpha \xi} & \xi \end{bmatrix}.
\end{equation*}
This representation covers every positive semidefinite
operator on $S_{\mat{V}}$ satisfying
$\langle\mU,\C_{\mat{V}}(\mU)\rangle=\xi$,
up to the sign of $\mat{V}$.
The additional condition that
$\frac12\B_{\mat{V}}-\C_{\mat{V}}$ is singular
determines the relationship between $\alpha$ and
$\sigma$. As shown below, the choice $\sigma=0$
then recovers the rank-one surrogate
$\C^\star_{\mat{V}}$.
We establish the desired comparison by showing
that, for each fixed $\omega\in[0,2]$,
the smallest eigenvalue of
$\mat{B}_{\mat{V}}-\omega\mat{C}_{\mat{V}}$
is nondecreasing in $\sigma$.

Since $\frac12\B_{\mat{V}}-\C_{\mat{V}}$ is
positive semidefinite and singular by assumption,
its determinant vanishes:
\begin{equation*}
    \det \left(\frac{1}{2}\mat{B}_{\mat{V}} - \mat{C}_{\mat{V}} \right) = \begin{vmatrix}
        \pi_{\mat{V}} - \alpha - \sigma & -\sqrt{\alpha \xi}\\ 
        -\sqrt{\alpha \xi} & \mu - \xi
        \end{vmatrix} = (\pi_{\mat{V}} - \alpha - \sigma) (\mu - \xi) - \alpha \xi = 0.
\end{equation*}
This poses the first constraint relating $\sigma$ and $\alpha$: 
\begin{equation}\label{eq:constraint1}
\alpha = (\pi_{\mat{V}} - \sigma) \left(1 - \frac{\xi}{\mu}\right).
\end{equation}
Moreover, since $\C_{\mat{V}}\preccurlyeq
\frac12\B_{\mat{V}}$, we have
$\alpha+\sigma\leq\pi_{\mat{V}}$,
and hence $0\leq\sigma\leq\pi_{\mat{V}}$.
Conversely, throughout this interval,
the constraint~(\ref{eq:constraint1})
defines a positive semidefinite
$\mat{C}_{\mat{V}}$ for which
$\frac12\mat{B}_{\mat{V}}-\mat{C}_{\mat{V}}$
remains positive semidefinite and singular.

For any relaxation value $\omega$, we consider the smaller eigenvalue of $\mat{B}_{\mat{V}} - \omega \mat{C}_{\mat{V}}$:
\begin{align*}
    \lambda_1(\mat{B}_{\mat{V}} - \omega \mat{C}_{\mat{V}}) = 
    \lambda_1\left(\begin{bmatrix} 2\pi_{\mat{V}} -\omega(\alpha + \sigma) & - \omega\sqrt{\alpha \xi}\\
    - \omega\sqrt{\alpha \xi} & 2\mu - \omega \xi\end{bmatrix} \right).
\end{align*}
Since eigenvalues of a $2 \times 2$ symmetric matrix $\mat{M}$ are $\left(\tr \mat{M}/2\right) \pm \sqrt{\left(\tr \mat{M}/2\right)^2 - \det \mat{M}}$, we introduce these quantities as:
\begin{align*}
    X &:= \tr{\left(\mat{B}_{\mat{V}} - \omega \mat{C}_{\mat{V}}\right)}/2 = \pi_{\mat{V}} + \mu - \frac{1}{2}\omega(\alpha + \sigma + \xi)\\
    Y &:= \det{\left(\mat{B}_{\mat{V}} - \omega \mat{C}_{\mat{V}}\right)} = 2\pi_{\mat{V}}\mu (2 - \omega) - \omega(2-\omega)\sigma\xi. %
\end{align*}
The last equality follows from the constraint~(\ref{eq:constraint1}). The smaller eigenvalue of $\mat{B}_{\mat{V}} - \omega \mat{C}_{\mat{V}}$ is given by $\lambda_1(\mat{B}_{\mat{V}} - \omega \mat{C}_{\mat{V}}) = X - \sqrt{X^2 - Y}$.

We first consider $0<\omega<2$, $\xi<\mu$, and
$0\leq\sigma<\pi_{\mat{V}}$. In this range,
$\alpha>0$ and $\xi>0$, so the off-diagonal
entry of $\mat{B}_{\mat{V}}-\omega\mat{C}_{\mat{V}}$
is nonzero. Consequently, its two eigenvalues
are distinct, and differentiation is justified.

The boundary cases are straightforward:
when $\omega=0$ or $\omega=2$, the smaller
eigenvalue equals $2\mu$ or $0$, respectively,
independently of $\sigma$. When $\xi=\mu$,
it equals $(2-\omega)\mu$, again independently
of $\sigma$. The endpoint
$\sigma=\pi_{\mat{V}}$ follows by continuity.

For the nondegenerate case, differentiation
with respect to $\sigma$ gives:
\begin{align*}
    \frac{\partial \lambda_1(\mat{B}_{\mat{V}} - \omega \mat{C}_{\mat{V}})}{\partial \sigma} = \frac{\partial X}{\partial \sigma} - \frac{1}{2\sqrt{X^2 - Y}}\left(2X \frac{\partial X}{\partial \sigma} - \frac{\partial Y}{\partial \sigma} \right).
\end{align*}
By using the constraint in~(\ref{eq:constraint1}), we have $\frac{\partial X}{\partial \sigma} = -\frac{\omega \xi}{2\mu}$. Similarly, $\frac{\partial Y}{\partial \sigma} = \omega (\omega - 2) \xi$ that yields:
\begin{align*}
    \frac{\partial \lambda_1(\mat{B}_{\mat{V}} - \omega \mat{C}_{\mat{V}})}{\partial \sigma} &= -\frac{\omega \xi}{2\mu} - \frac{1}{2\sqrt{X^2 - Y}}\left(-2X\frac{\omega \xi}{2\mu} - \omega(\omega - 2) \xi \right)\\ &= \frac{\omega\xi}{2\mu}\left(\frac{X-(2-\omega)\mu}{\sqrt{X^2 - Y}} - 1 \right).
\end{align*}
By our assumption that
$\frac12\B_{\mat{V}}-\C_{\mat{V}}$ is
positive semidefinite, we have
$\C_{\mat{V}}\preccurlyeq\frac12\B_{\mat{V}}$
and   $\lambda_1(\mat{B}_{\mat{V}} - \omega \mat{C}_{\mat{V}}) \ge \lambda_1(\mat{B}_{\mat{V}} - \omega\frac12\mat{B}_{\mat{V}}) = (2 - \omega) \mu$. Since $X = \tr{\left(\mat{B}_{\mat{V}} - \omega \mat{C}_{\mat{V}}\right)}/2$ is the average of its two eigenvalues and $(2-\omega)\mu$ is a lower bound for the smaller eigenvalue, we have $X - (2 - \omega)\mu \ge 0$. To show $X - (2 - \omega)\mu \ge \sqrt{X^2 - Y}$ we take the square and consider the difference:
\begin{align*}
    (X - (2 - \omega)\mu)^2 - X^2 + Y &= (2-\omega)^2 \mu^2 - 2(2 - \omega)\mu X + Y \\
    &= \omega(2-\omega)\mu \left[\alpha-(\mu - \xi) +\sigma \left(1 - \frac{\xi}{\mu}\right) \right]\\
    &= \omega(2-\omega)\mu \left[(\pi_{\mat{V}} - \sigma) \left(1 - \frac{\xi}{\mu}\right)-(\mu - \xi) +\sigma \left(1 - \frac{\xi}{\mu}\right) \right]\\
    &= \omega(2-\omega)\mu \left[\pi_{\mat{V}}\left(1 - \frac{\xi}{\mu}\right)-(\mu - \xi) \right]\\
    &= \omega(2 - \omega)(\mu - \xi)(\pi_{\mat{V}} - \mu) \ge 0.
\end{align*}
Lemma~\ref{lem:ingredients} shows $\xi \le \mu$ and $\mu \le \pi_{\mat{V}}$ since $2\mu = \langle \mU, \B(\mU)\rangle \le \langle \mat{V}, \B(\mat{V})\rangle = 2\pi_{\mat{V}}$ as $\mat{V}$ is orthogonal to $\mU$. 
Consequently, the smallest eigenvalue is
nondecreasing in $\sigma$ for every
$\omega\in[0,2]$. Since $\sigma=0$ corresponds
to the rank-one surrogate and the saturated
compression $\widehat{\C}_{\mat{V}}$ corresponds
to some $\sigma\geq0$, we obtain
\[
\lambda_{\min}(\B_{\mat{V}}
-\omega\C^\star_{\mat{V}})
\leq
\lambda_{\min}(\B_{\mat{V}}
-\omega\widehat{\C}_{\mat{V}})
\leq
\lambda_{\min}(\B_{\mat{V}}
-\omega\C_{S_{\mat{V}}}).
\]
The second inequality follows from
$\widehat{\C}_{\mat{V}}\succeq
\C_{S_{\mat{V}}}$ and $\omega\geq0$.
This proves the proposition.
\end{proof}

\begin{proof}[Proof of Proposition~\ref{pro:c_star}]
    We recall the $2 \times 2$ matrix representations $\mat{B}_{\mat{V}}$ and $\mat{C}^\star_{\mat{V}}$ from~(\ref{eq:c_star_v}) and denote $\pi_{\mat{V}} := \frac{1}{2}\langle \mat{V}, \B(\mat{V}) \rangle$. 
By the min-max theorem,
\[
\min_{\substack{\|\mat{V}\|_F=1\\
\mat{V}\perp\mU}}\pi_{\mat{V}}
=\frac12\lambda_2(\B),
\]
and this minimum is attained at $\mat{V}=\mU^\prime$.
 In particular, choosing $\mat{V}=\mU'$ recovers $\C^\star_{\mat{V}}=\C^\star$. 
    We establish the desired comparison by showing that,
    for each fixed $\omega\in[0,2]$, the smallest
    eigenvalue of
    $\mat{B}_{\mat{V}}-\omega\mat{C}^\star_{\mat{V}}$
    is nondecreasing as $\pi_{\mat{V}}$ increases from its
    minimum value $\frac12\lambda_2(\B)$.
    As in the preceding proof, we establish this
    by showing that the derivative with respect
    to $\pi_{\mat{V}}$ is nonnegative.
     The smaller eigenvalue of $\mat{B}_{\mat{V}} - \omega \mat{C}^\star_{\mat{V}}$ as a function of $\pi_{\mat{V}}$ is:
\[
\lambda_1(\mat{B}_{\mat{V}}-\omega\mat{C}^\star_{\mat{V}})
=
\lambda_1\left(
\begin{bmatrix}
2\pi_{\mat{V}}-\omega\gamma_{\mat{V}}
& -\omega\sqrt{\gamma_{\mat{V}}\xi}\\
-\omega\sqrt{\gamma_{\mat{V}}\xi}
& 2\mu-\omega\xi
\end{bmatrix}
\right).
\]
recalling from~(\ref{eq:c_star_v}) that $\gamma_{\mat{V}} = \pi_{\mat{V}} \left(1 - \xi/\mu\right)$. As before, eigenvalues of this $2 \times 2$ symmetric matrix can be derived from its trace and determinant:
\begin{align*}
    X &:= \frac12\tr{\left(\mat{B}_{\mat{V}} - \omega \mat{C}^\star_{\mat{V}}\right)} = \pi_{\mat{V}} + \mu - \frac{\omega}{2}(\gamma_{\mat{V}} + \xi) = \pi_{\mat{V}} \left(1 - \frac{\omega}{2}\left(1 -\frac{\xi}{\mu}\right) \right) + \mu - \frac{\omega \xi}{2}\\
    Y &:= \det{\left(\mat{B}_{\mat{V}} - \omega \mat{C}^\star_{\mat{V}}\right)} = 2\pi_{\mat{V}}\mu(2 - \omega).
\end{align*}
The smaller eigenvalue of
$\mat{B}_{\mat{V}}-\omega\mat{C}^\star_{\mat{V}}$
is $X-\sqrt{X^2-Y}$.

We first consider the boundary cases.
When $\omega=0$, the smaller eigenvalue is
$2\mu$, independently of $\pi_{\mat{V}}$.
When $\omega=2$, it is zero, since
$\frac12\B_{\mat{V}}-\C^\star_{\mat{V}}$
is positive semidefinite and singular.
If $\xi=\mu$, then $\gamma_{\mat{V}}=0$,
and the smaller eigenvalue is
$(2-\omega)\mu$, again independently
of $\pi_{\mat{V}}$.

It therefore suffices to consider
$0<\omega<2$ and $\xi<\mu$.
Since $\pi_{\mat{V}}\geq\mu>0$ and
$\mu^2\leq\xi<\mu$, we have
$\gamma_{\mat{V}}
=\pi_{\mat{V}}(1-\xi/\mu)>0$
and $\xi>0$. Thus, the off-diagonal
entry is nonzero,
so the two eigenvalues are distinct
and differentiation is justified.
Writing
$Z:=\partial X/\partial\pi_{\mat{V}}
=1-\frac{\omega}{2}(1-\xi/\mu)$
and using
$\partial Y/\partial\pi_{\mat{V}}
=2\mu(2-\omega)$, we obtain:
\begin{align*}
    \frac{\partial \lambda_1(\mat{B}_{\mat{V}} - \omega \mat{C}^\star_{\mat{V}})}{\partial \pi_{\mat{V}}} = Z - \frac{X Z - \mu (2 - \omega)}{\sqrt{X^2 - Y}}.
\end{align*}
We need to prove that $\partial \lambda_1(\mat{B}_{\mat{V}} - \omega \mat{C}^\star_{\mat{V}})/{\partial \pi_{\mat{V}}} \ge 0$, alternatively: $Z\sqrt{X^2 - Y} \ge XZ - \mu(2 - \omega)$. Note that from $\omega \in [0, 2]$ and Lemma~\ref{lem:ingredients}, we have $Z \ge 0$ and if $X Z - \mu(2 - \omega) \le 0$ the positivity is trivially established. With the assumption that the latter is positive, we square both sides:
\begin{align*}
    Z^2 (X^2 - Y) &\ge X^2Z^2 + \mu^2(2 - \omega)^2 - 2\mu(2 - \omega) X Z \quad {\rm or}\\
    2\mu(2 - \omega) X Z &\ge \mu^2 (2 - \omega)^2 + Z^2 Y.
\end{align*}
Plugging in the value of $Y$, we need to show:
\begin{align*}
    2 X Z \ge \mu (2 - \omega) + 2 \pi_{\mat{V}} Z^2.
\end{align*}
Noting that $X = \pi_{\mat{V}} Z + \mu - \omega\xi/2$, the inequality becomes:
\begin{align*}
    2\mu Z - \omega \xi Z \ge \mu (2 - \omega).
\end{align*}
Bringing in $Z$ results in the inequality:
\begin{align*}
    \frac{\omega^2 \xi}{2} \left(1 - \frac{\xi}{\mu}\right) \ge 0,
\end{align*}
which holds because of Lemma~\ref{lem:ingredients}.

Therefore, the smallest eigenvalue of
$\mat{B}_{\mat{V}}-\omega\mat{C}^\star_{\mat{V}}$
is nondecreasing in $\pi_{\mat{V}}$.
Since $\pi_{\mat{V}}$ is minimized by
$\mat{V}=\mU^\prime$, we conclude that
\[
\lambda_{\min}(\B_{S^\star}
-\omega\C^\star_{S^\star})
\leq
\lambda_{\min}(\B_{\mat{V}}
-\omega\C^\star_{\mat{V}})
\]
for every $\omega\in[0,2]$.
This proves the proposition.
\end{proof}
\bibliographystyle{amsplain}
\bibliography{references}
\end{document}